\documentclass[twoside, 11pt]{article}
\usepackage{indentfirst,bm,latexsym,amsfonts,amssymb,amsmath,mathrsfs}
\RequirePackage{ifthen,calc}
\usepackage{theorem}
\usepackage{color}

 \catcode`@=11
 \def\@evenhead{\hbox to\textwidth{\footnotesize\rm\thepage \hfill
  {\it }}} 

 \def\@oddhead{\hbox to \textwidth{\footnotesize{\it
  On the barrier problem of branching random walk in a
time-inhomogeneous random environment  } \hfill\thepage}}

\catcode`@=12
\newcommand\ack{\section*{Acknowledgement}}

\newtheorem{thm}{\noindent Theorem}[section]
\newtheorem{lem}[thm]{\noindent Lemma}
\newtheorem{cor}[thm]{\noindent Corollary}
\newtheorem{prop}[thm]{\noindent Proposition}

\newtheorem{remark}[thm]{\noindent Remark}
\newtheorem{lemma}[thm]{Lemma}
\newtheorem{example}[thm]{Example}
\theoremheaderfont{\normalfont\bfseries} \theorembodyfont{\slshape}
{\theorembodyfont{\rmfamily}}
{\theorembodyfont{\rmfamily}\newtheorem{defn}{\noindent Definition}[section]}
{\theorembodyfont{\rmfamily}}
{\theorembodyfont{\rmfamily}}
{\theorembodyfont{\rmfamily}\newtheorem{rem}{\noindent Remark}[section]}
{\theorembodyfont{\rmfamily}}
{\theorembodyfont{\rmfamily}}
\usepackage[colorlinks,citecolor=blue,linkcolor=blue]{hyperref}

 \def\beqlb{\begin{eqnarray}}\def\eeqlb{\end{eqnarray}}
 \def\beqnn{\begin{eqnarray*}}\def\eeqnn{\end{eqnarray*}}
 \newcommand{\bgeqn}{\begin{equation}}
\newcommand{\edeqn}{\end{equation}}
\def\ra{\rightarrow}

 \numberwithin{equation}{section}
\def\qed{\hfill$\square$\smallskip}
\def\no{\nonumber}
\def\L{{\mathcal L}}

\def\ra{\rightarrow}

\def\bfE{{\mathbb{E}}}

\def\mbfE{{\mathbf{E}}}
\def\mbr{{\mathbb{R}}}
\def\bfP{{\mathbb{P}}}

\def\mbfP{{\mathbf{P}}}
\def\bfR{{\mathbb{R}}}

\def\bfN{{\mathbb{N}}}
\def\1{{\mathbf{1}}}

\def\ee{\varepsilon}

\begin{document}
\title{\bf The extinction rate of a branching random walk with a barrier in a time-inhomogeneous random environment
}
\author{You Lv\thanks{Email: lvyou@dhu.edu.cn }~~
\\
\\ \small College of Science, Donghua University,
\\ \small Shanghai 201620, P. R. China.
}
\date{}
\maketitle


\noindent\textbf{Abstract}: Consider a supercritical branching random walk in a time-inhomogeneous random environment. We impose a selection (called barrier) on survival in the following way. The position of the barrier may depend on the generation and the environment. 
In each generation, only the individuals born below the barrier can survive and reproduce. When the barrier causes the extinction of the system, we give the extinction rate in the sense of $L^p~(p\geq 1)$. Moreover, we show the $L^p$ convergence of the small deviation probability for a random walk with random environment in time.

\smallskip


\smallskip


\smallskip

\section{Introduction and result}
\subsection{Model}
We consider a branching random walk on $\mathbb{R}$ in a time-inhomogeneous i.i.d. random environment
(BRWre), which is an extension of the time-homogeneous branching random walk (BRW). For a BRW, the reproduction law (including displacement and branching) of each generation is determined by a common point process, while for a BRWre, the reproduction law of each generation is sampled independently according to a common distribution on the collection of the point processes on $\bfR.$ The mathematical definition is as follows. 

Let $(\Pi,\mathcal{F}_{\Pi})$ be a measurable space and $\Pi\subseteq\tilde{\Pi}:=\{\mathfrak{m}:\mathfrak{m}~\text{is~a~point~process~on~}\bfR\}.$ The random environment $\L$ is defined as an i.i.d. sequence of random
 variables $\{\L_1$,~$\L_2$,~$\cdots,\L_n,\cdots\}$, where $\L_1$ takes values in $\Pi$.
 Let $\nu$ be the law of $\L$, then we call the product space $(\Pi^{\bfN}, \mathcal{F}_{\Pi}^{^{\bigotimes}\bfN}, \nu)$ the \emph{environment space}. For any realization $L:=\{L_1$,~$L_2$,~$\cdots,L_n,\cdots\}$ of $\L$, a time-inhomogeneous branching random walk driven by the environment $L$ is a process constructed as follows.

(1)~At time $0,$ an initial particle $\phi$ in generation $0$ is located at the origin.

(2)~At time $1,$ the initial particle $\phi$ dies and gives birth to $N(\phi)$ children who form the first generation. These children are located at $\zeta_i(\phi), 1\leq i\leq N(\phi),$ where the distribution of the random vector $X(\phi):=(N(\phi), \zeta_1(\phi),\zeta_2(\phi),\ldots)$ is $L_1.$ 

(3)~ Similarly, at generation $n+1,$ every particle $u$ alive at generation $n$ dies and gives birth to $N(u)$ children. If we denote $\zeta_i(u), 1\leq i\leq N(u)$ the displacement of the children with respect to their parent $u$, then the distribution of $X(u):=(N(u), \zeta_1(u),\zeta_2(u),\cdots)$ is $L_{n+1}.$ We should emphasize that conditionally on any given environment $L,$ all particles in this system reproduce independently. 

Conditionally on $\L,$ we write $(\Gamma,\mathcal{F}_{\Gamma}, \mbfP_{\L})$ for the probability space under which the
time-inhomogeneous branching random walk is defined. The probability $\mbfP_{\L}$ is usually called a {\it quenched law}.  We define the probability $\mathbf{P}:=\nu\bigotimes\mbfP_{\L}$ on the product space $(\Pi^{\bfN}\times\Gamma,\mathcal{F}_{\Pi}^{^{\bigotimes}\bfN}\bigotimes\mathcal{F}_{\Gamma})$ such that
\begin{eqnarray}\label{APP}\mathbf{P}(F\times G)=\int_{\L\in F}\mbfP_{\L}(G)~d\nu(\L), ~~F\in \mathcal{F}_{\Pi}^{^{\bigotimes}\bfN}, ~G\in\mathcal{F}_{\Gamma}.\end{eqnarray}
The marginal distribution of probability $\mathbf{P}$ on $\Gamma$ is usually called an {\it annealed law}. 
Throughout this paper, we consider the case $F=\Pi^{\bfN}.$ Hence without confusion we also denote $\mathbf{P}$ the annealed law and abbreviate $\mathbf{P}(\Pi^{\bfN}\times G)$ to~$\mathbf{P}(G).$ Moreover, we write $\mbfE_{\L}$ and~$\mathbf{E}$ for the corresponding expectation of $\mbfP_{\L}$ and~$\mathbf{P}$ respectively.

We write $\mathbf{T}$ for the (random) genealogical tree of the process. For a given particle $u\in\mathbf{T}$ we write $V(u)\in\bfR$ for the position of $u$ and $|u|$ for the generation at which $u$ is alive. Then $(\mathbf{T}, V, \mbfP_{\L}, \mathbf{P})$ is called the {\it branching random walk in the time-inhomogeneous random environment $\L$} (BRWre). This model was first introduced in Biggins and Kyprianou \cite{BK2004}. If there exists a point process $\iota\in\Pi$ such that $\mathbf{P}(\L_1=\iota)=1$ (thus $\mathbf{P}(\L_i=\iota)=1, \forall i\in\bfN^+:=\{1,2,\cdots,n,\cdots\}$), then we usually call the environment a {\it degenerate environment} 
 and the BRWre degenerates to a BRW.

Denote $M_n:=\min_{|u|=n}V(u)$ the minimal displacement in generation $n$ and $$\kappa_n(\theta):=\log\mbfE_{\L}\left(\sum^{N(u)}_{i=1}e^{-\theta \zeta_i(u)}\right),~|u|=n-1,~ \theta\in[0,+\infty)$$ the log-Laplace transform function of the random point process $\L_n.$
Throughout the present paper, we assume that
\beqlb\label{sect-2}\exists\theta_*>0, ~~\kappa(\theta_*)<+\infty,~~\exists \vartheta\in(0,\theta_*),~~\kappa(\vartheta)=\vartheta\kappa'(\vartheta), ~~\kappa(0)\in(0, +\infty),\eeqlb
where $\kappa(\theta):=\mbfE(\kappa_n(\theta))$\footnote{Note that \eqref{sect-2} ensures that $\kappa(\cdot)$ is well-defined in $[0,\theta_*].$ Moreover, for any $i,j\in\bfN, \theta>0, \mbfE(\kappa_i(\theta))=\mbfE(\kappa_j(\theta))$ since $\L_1,...\L_n,...$ is an i.i.d. sequence.} and $\kappa'$ is the derivative of $\kappa.$ \eqref{sect-2} is also a basic assumption in the following papers which have contributed some deep results in BRWre.
Huang and Liu \cite{HL2014} proved
\beqlb\label{HL}\lim\limits_{n\rightarrow\infty}\frac{M_n}{n}=-\frac{\kappa(\vartheta)}{\vartheta},~~~ \mathbf{P}-{\rm a.s.},\eeqlb
and obtained large deviation principles for the counting measure of the BRWre. Conclusions on the central limit theorem of the BRWre can be found in Gao et al. \cite{GLW2014} and Gao, Liu \cite{GL2016}. A moderate deviation principle for the counting measure was investigated in Wang, Huang \cite{WH2017}.~~

The second order of the asymptotic behavior of $M_n$ was given in Mallein, Mi{\l}o\'{s} \cite{MM2016}. They showed that  \beqlb\label{SOBRWre}\lim\limits_{n\rightarrow\infty}\frac{M_n+\vartheta^{-1} K_n}{\log n}=c,~~~\text{in~~Probability}~~ \mathbf{P},\eeqlb
where $c$ is an explicit constant and \beqlb\label{Kn}K_n:=\sum_{i=1}^n\kappa_i(\vartheta),~~~~ K_0:=0.\eeqlb
\eqref{SOBRWre} provides a basis for the research on the barrier problem, see the footnote in the next subsection.

\subsection{Barrier problem}
In the present paper, we focus on a barrier problem of BRWre. 
The motivation to consider the barrier problem is from the research on parallel simulations in Lubachevsky et al. \cite{LSW1989,LSW1990}. The BRW with barrier was first introduced in Biggins et al. \cite{BLSW1991}. The so-called ``barrier" is a function $\varphi:\bfN\ra\bfR.$ For any particle $u\in\mathbf{T},$ we will erase it and all its descendants as long as $V(u)>\varphi(|u|).$ The new branching particle system after removing is called a BRW with a barrier $\varphi.$
Assume the underlying Galton-Watson process is supercritical.
The barrier problem of BRW was raised mainly in the following two aspects.

      1. Consider the impact (extinction/survival) of the barrier on the particle system, see \cite{BLSW1991,GHS2011,BJ2012}.

      2. Consider the extinction rate when the survival probability is 0, see \cite{AJ2011,BJ2012}.

      3. Consider the total progeny when the survival probability is 0, see \cite{A2010,AB2011}.

In the present paper we study the aforementioned second point for BRWre with a random barrier in the sense of $L^p (p\geq 1).$

Let us first introduce some notation for a better expression of the barrier problem of BRWre. On the tree $\mathbf{T}$ we define a partial order $>$ such that $u>v$ if $v$ is an ancestor of $u$. We write $u\geq v$ if $u>v$ or $u=v$ (i.e., $u$ and $v$ are the same one). We define an {\it infinite path} $u_{\infty}$ through $\mathbf{T}$ as a sequence of particles $u_{\infty}:=(u_i,i\in\bfN)$ such that $$\forall i\in\bfN,~~ |u_i|=i,~~ u_{i+1}>u_{i},~~ u_0=\phi~(\text{the initial particle}).$$ For any $i\leq |u|,$ we usually write $u_i$ for the ancestor of $u$ in generation $i.$ Let $\mathbf{T}_n:=\{u\in\mathbf{T}:|u|=n\}$ be the set of particles in generation $n$ and $\mathbf{T}_\infty$ the collection of all infinite paths through $\mathbf{T}.$
Then we see that the event $$\mathcal{S}:=\{\exists u_{\infty}:=(u_0,u_1,u_2, \ldots u_n, \ldots)\in \mathbf{T_{\infty}}, \forall i\in\bfN,  V(u_i)\leq \varphi_{_\L}(i)\}$$ represents that the system still survives after we add the barrier $\varphi_{_\L},$ where we add a subscript $\L$ to $\varphi$ since the barrier we consider may depend on the random environment. Denote
$$Y_n:=\sharp\{|u|=n:~\forall i\leq n, ~V(u_i)\leq \varphi_{\L}(i)\}$$
the size of surviving population at generation $n.$

 By the light\footnote{From the definition of the barrier, it is reasonable to image that the quenched survival probability $\mbfP_{\L}(\mathcal{S})$ will approach $0$ when the barrier is close to the trajectory of $M_n.$ Comparing \eqref{HL} with \eqref{SOBRWre}, we see that $M_n$ is closer to $\vartheta^{-1}K_n$ than $\vartheta^{-1}n\mbfE(K_1)$ for $n$ large enough. Therefore, we set a random barrier rather than a constant one in the context of the BRW---for a BRW with a barrier, the barrier in generation $n$ is usually set as $d^*n+dn^{\alpha}$, where $d^*$ is the limit of $M_n/n$, see \cite{AJ2011,BLSW1991,GHS2011,BJ2012}.} of  
$\eqref{SOBRWre}$, Lv, Hong \cite{LY3} added a barrier function $\varphi_{\L}(i):=-\vartheta^{-1} K_i+di^{\alpha}$ to the BRWre and obtained the following result.
\begin{thm}\label{as} (Lv $\&$ Hong \cite[Theorem 2.6 (2a) and (2b)]{LY3})~If \eqref{sect-2} holds and there exist $\lambda_0>6,\lambda_1>3,\lambda_{2}>2,\lambda_{3}>6,\lambda_{4}>0, \lambda_{5}\leq -1, \lambda_6\geq 1$ such that
\begin{eqnarray}\label{c2'a}
\mathbf{E}\left(|\kappa_1(\vartheta)-\vartheta \kappa'_1(\vartheta)|^{\lambda_0}\right)<+\infty;
\end{eqnarray}
\begin{eqnarray}\label{c2'b}
\mbfE\left(\left[\frac{\mbfE_{\L}\left(\sum_{i=1}^{N(\phi)}|\zeta_i(\phi)+\kappa'_1(\vartheta)|^{\lambda_2}e^{-\vartheta \zeta_i(\phi)}\right)}{\mbfE_{\L}\left(\sum_{i=1}^{N(\phi)}e^{-\vartheta \zeta_i(\phi)}\right)}\right]^{\lambda_1}\right)<+\infty.
\end{eqnarray}
\begin{eqnarray}\label{c3a}\mbfE(|\kappa_1(\vartheta+\lambda_4)|^{\lambda_3})+\mbfE(|\kappa_1(\vartheta)|^{\lambda_3})<+\infty,~~ \mbfE([\log^+\mbfE_{\L}(N(\phi) ^{1+\lambda_4})]^{\lambda_3})<+\infty,\end{eqnarray}
where $\log^+\cdot:=\log\max(1,\cdot), ~\log^-\cdot:=|\log\min(1,\cdot)|,$ and

\begin{eqnarray}\label{T<}
\mathbf{E}\left(\left[\log^-\mbfE_{\L}\left(\1_{N(\phi)\leq |\lambda_5|}\sum_{i=1}^{N(\phi)}\1_{\{\vartheta \zeta_i(\phi)+ \kappa_1(\vartheta)\in [\lambda_5, \lambda^{-1}_5]\}} \right)\right]^{\lambda_6}\right)<+\infty,\end{eqnarray}
 then we have the extinction rates as follows. 

{\rm (1)}~If $\alpha=\frac{1}{3}, d\in(0,d_c),$ then there exists a negative constant $b_1$ depending on $d$ such that \begin{eqnarray}\label{b1}
\lim\limits_{n\rightarrow\infty}\frac{\log\mbfP_{\L}(Y_n>0)}{\sqrt[3]{n}}=b_1,
\end{eqnarray}
holds in the sense of~${\rm \mathbf{P}-a.s.}$, where $d_c$ is a positive constant. 

{\rm (2)}~If~$\alpha\in(0,\frac{1}{3}), d\geq 0$, then there exists a negative constant $b_2$ (not depending on $d$) such that \begin{eqnarray}\label{b2}
\lim\limits_{n\rightarrow\infty}\frac{\log\mbfP_{\L}(Y_n>0)}{\sqrt[3]{n}}=b_2, 
\end{eqnarray}
holds in the sense of~${\rm \mathbf{P}-a.s.}$

\end{thm}

The explicit expressions of $d_c, b_1$ and $b_2$ in the above theorem has be obtained in \cite{LY3}. Here we do not give the expressions since they are not involved in the present paper. 
\cite{LY3} proved that under the assumptions \eqref{sect-2}, \eqref{c2'a}-\eqref{T<} with $\lambda_6>2$, $\mbfP_{\L}(\mathcal{S})>0,~{\rm \mathbf{P}-a.s.}$ when $\alpha>\frac{1}{3}, d>0$ or $\alpha=\frac{1}{3}, d>d_c$; under the assumptions \eqref{sect-2}, \eqref{c2'a} and \eqref{c2'b}, $\mbfP_{\L}(\mathcal{S})=0,~{\rm \mathbf{P}-a.s.}$ when $\alpha=\frac{1}{3}, d<d_c$ or $\alpha<\frac{1}{3}.$
In fact, \cite[Theorem 2.5 and 2.6]{LY3} extended the main results in \cite{AJ2011} and \cite{BJ2012} to the case of random environment. We also refer to \cite{LY3} for some detailed explanations on conditions  \eqref{sect-2} and \eqref{c2'a}-\eqref{T<} and an example satisfying all the conditions.

The present paper looks for the sufficient conditions for the $L^p$ convergence in \eqref{b1} and \eqref{b2}. As a basis for $L^p$ convergence, we first give a group of sufficient conditions for convergence in probability.


\begin{thm}\label{inp}
If \eqref{sect-2}, \eqref{c2'a}-\eqref{T<} hold with constants $\lambda_0>3,\lambda_1>2,\lambda_{2}>2,\lambda_{3}>3,\lambda_{4}>0, \lambda_{5}\leq -1, \lambda_{6}\geq 1,$
then the convergence in \eqref{b1} and \eqref{b2} hold in probability $\mathbf{P}.$
\end{thm}

  We omit the proof of this theorem since it is very similar to the proof of Theorem \ref{as}. In \cite{LY3}, we show \eqref{b1} and \eqref{b2} by many-to-one formula \cite[Lemma 3.1]{LY3} (transfer the BRWre to a random walk with random environment in time, abbreviated as RWre, see Section 4 for the definition) and the small deviation principle \cite[Theorem 2]{Lv201802} (see also Theorem \ref{mog}(b) in the present paper) for RWre.  In Lv, Hong \cite{Lv201802} we gave two groups of sufficient conditions for the limit behavior of the scaling small deviation probability in the sense of almost surely and in probability respectively. For the proof of Theorem \ref{inp}, 
the method used in \cite[Theorem 2.6]{LY3} still works as long as we replace  the ``almost surely" version of the small deviation principle (Theorem \ref{mog}(b)) by the ``in probability" version (Theorem \ref{mog}(a)).

\subsection{Main result}

  The following theorem is the main result in the present paper.

\begin{thm}\label{Lp}
Assume that \eqref{sect-2}, \eqref{c2'a}-\eqref{T<} hold with constants
\begin{eqnarray}\label{assu1}\lambda_0>3,~~\lambda_1>2,~~\lambda_{2}>2,~~\lambda_{3}>3,~~\lambda_{4}>0,~~\lambda_{5}\leq-1,~~ \lambda_6\geq 1,~~\frac{\sigma^2}{\sigma^2_*}< \frac{\lambda_2-2}{\lambda_0-2}\end{eqnarray}
and
\begin{eqnarray}\label{T>}
\mathbf{E}\left(\left[\log^-\mbfE_{\L}\left(\1_{N(\phi)\leq |\lambda_5|}\sum_{i=1}^{N(\phi)}\1_{\{\vartheta \zeta_i(\phi)+ \kappa_1(\vartheta)\in [0, |\lambda_5|]\}} \right)\right]^{\lambda_6}\right)<+\infty,\end{eqnarray}
where $$\sigma^2:=\mathbf{E}\left(\Big(\kappa_1(\vartheta)-\vartheta\kappa'_1(\vartheta)\Big)^2\right),~~\sigma^2_*:=\vartheta^2\mathbf{E}(\kappa''_1(\vartheta)).$$ 
If constants $p, t$ satisfy that
 \begin{eqnarray}\label{assu2}p\geq 1,~t\geq 1,~t\lambda_6\geq 2,~ p\in\left[\lambda_6-\frac{1}{t},\lambda_6\right),~\min\left\{\frac{\lambda_0}{2},\lambda_1,\frac{\lambda_3}{2}\right\}>\frac{\lambda_6}{\lambda_6-p},\end{eqnarray}
then the convergence in \eqref{b1} and \eqref{b2} hold in $L^p.$
\end{thm}

The $L^p$ convergence reflects the characteristic of studying random environment.  Note that for BRW, the convergence in \eqref{b1} and \eqref{b2} is essentially the convergence of a sequence of numbers since the probability of $\{Y_n>0\}$ has no randomness. But for BRWre, the quenched probability $\mbfP_{\L}(Y_n>0)$ is a conditional probability and the law of $\mbfP_{\L}(Y_n>0)$ will be totally determined by the random environment hence both almost surely convergence and $L^p$ convergence for \eqref{b1} and \eqref{b2} are of independent interest and challenging. In \cite{LY3}, the  almost surely convergence has been obtained as what we have mentioned in Theorem \ref{as}. In the present paper we focus on the $L^p$ convergence. It should be noted that the situation and difficulty of $L^p$ convergence will be quite different from that of the almost surely convergence.

Let us give a brief explanation about the difference. For the almost surely convergence, we need to show the probability of the extreme environment (in other words, the bad environment) is small enough. But this is not enough for $L^p$ convergence since the impact of the extreme environment may be also extremely awful. Therefore, we need to pick out the extreme environment accurately and estimate the influence of the extreme environment carefully. The difficulty to show Theorem \ref{Lp} mainly focuses on how to estimate the ($k$-th) moment of a kind of $\log$-quenched probability, see \eqref{goal} for details.~Overcoming the difficulty is not only meaningful in the proof of Theorem \ref{Lp}, but also helpful in studying another barrier problem, see the two paragraphs below Remark \ref{rema3} for details\footnote{We postpone the details to Section 4 because some necessary notation for a clear explanation has not been introduced until Section 4.}.


Let us give some comments on the assumptions in Theorem \ref{Lp} as the ending of this subsection. First, \eqref{sect-2} is a basic assumption for BRW and varieties of its generalizations, which ensures the existence of the first order of $\min_{|u|=n}V(u)$ and the supercritical property of the underlying branching process. Assumptions \eqref{c2'a}, \eqref{c2'b}\footnote{
 \cite[Proposition 2.2]{LY3} showed that if there exists $\lambda>\frac{3}{2}$ such that
 ~$\mbfE((\kappa^{(4)}_1(\vartheta)+3[\kappa''_1(\vartheta)]^2)^{\lambda})<+\infty$,  where $\kappa^{(4)}_1(\vartheta):=\frac{d^4\kappa(\theta)}{d\theta}|_{\theta=\vartheta},$ then $\eqref{c2'b}$ holds.}, \eqref{c3a}, \eqref{assu1} and \eqref{assu2} all express that some mild integrability conditions on $\kappa_1(\vartheta), \kappa'_1(\vartheta)$ and $N(\phi)$ are required. Recalling the notation $\log^{-}\cdot:=|\log\min(\cdot,1)|,$ we see that the larger value $|\lambda_5|$ takes, the easier \eqref{T<} and \eqref{T>} hold. The following proposition and remark provide an intuitive rationalization for the assumptions \eqref{T<} and \eqref{T>}.
\begin{prop}\label{rem1+}
If the random environment is degenerate, both \eqref{T<} and \eqref{T>} can be deduced from \eqref{sect-2}.
\end{prop}

Since the barrier problem of BRW had always been considered under the corresponding assumptions of \eqref{sect-2} in the time-homogeneous case (e.g., \cite{AJ2011}, \cite{GHS2011}, \cite{BJ2012}),
this proposition explains that why these papers never set an assumption like \eqref{T<} or \eqref{T>}.

We prove this proposition after we introduce the many-to-one formula (see Section 2.1).
\begin{rem}\label{rema2}
\eqref{T<} can be roughly understood as that there exists a negative $\lambda_5$ near $-\infty$ such that
the distribution of $\mbfE_{\L}\left(\1_{N(\phi)\leq |\lambda_5|}\sum_{i=1}^{N(\phi)}\1_{\{\vartheta \zeta_i(\phi)+ \kappa_1(\vartheta)\in [\lambda_5, \lambda^{-1}_5]\}} \right)$ can not be too concentrated at a neighborhood of $0$. Let us consider an extreme circumstance that \begin{eqnarray}\label{qqq}\exists \ee>0 {\rm~such~that~} q:=\mbfP\left(\mbfE_{\L}\left(\sum_{i=1}^{N(\phi)}\1_{\{\vartheta \zeta_i(\phi)+ \kappa_1(\vartheta)\leq \ee\}} \right)=0\right)>0.\end{eqnarray}Obviously, \eqref{qqq} means that \eqref{T<} does not hold. Now we explain that \eqref{qqq} contradicts Theorem \ref{Lp}.
We should note that \eqref{qqq} is equivalent to saying that $\forall k\in\bfN^+, \mbfP(\mathcal{A}_k)>0$, where
$$\mathcal{A}_k:=\left\{\L: \mbfP_{\L}\left(\vartheta \min_{i\leq N(u)}\zeta_i(u)+\kappa_k(\vartheta)>\ee\right)=1, |u|=k-1\right\}.$$
We remind that $\mbfP(\cap_{k=1}^{n}\mathcal{A}_k)=\mbfP(\mathcal{A}_1)^n=q^n$ as the random environment is i.i.d.
Note that for $n$ large enough, $\{\L\in\cap_{k=1}^{n}\mathcal{A}_k\}\subset\{\min_{|u|=n}V(u)>-\vartheta^{-1}K_n+\ee n\}$, which means that
$$\mbfP(\mbfP_{\L}(Y_n>0)=0)\geq\mbfP(\mathcal{A}_1)^n>0.$$ That is to say, for any $n, \mbfE\left(n^{-1/3}\log\mbfP_{\L}(Y_n>0)\right)=+\infty,$ which means that the $L^p (p\geq 1)$ convergence in Theorem \ref{Lp} is not true.

Especially, if
 \begin{eqnarray}\label{nece}
\sharp\Pi<+\infty, ~~~\exists\varepsilon>0, ~~\mbfP\left(\mbfP_{\L}\left(\min_{i\leq N(\phi)}|\zeta_i(\phi)+\vartheta^{-1}\kappa_1(\vartheta)|>\varepsilon\right)=1\right)=1,\end{eqnarray} where $\sharp\cdot$ represents the number of elements in set $\cdot,$ then letting $\lambda_5\ra -\infty$, the monotone convergence theorem tells that the negative proposition of \eqref{T<} is equivalent to \eqref{qqq}. Hence under \eqref{nece}, we see \eqref{T<} is a necessary condition for Theorem \ref{Lp}.
\end{rem}

The examples in the coming subsection will lead us to see how these assumptions dominate our construction of the model intuitively.

\subsection{Examples}

Recall that \eqref{T<} and \eqref{qqq} are mutually exclusive. Now we give two examples: the first example satisfies \eqref{qqq} and  all the assumptions in Theorem \ref{Lp} except \eqref{T<}; the second example satisfies all the assumptions in Theorem \ref{Lp}.

We remind that we have given an example in \cite[Section 2]{LY3} which satisfies all the assumptions in Theorem \ref{as}. Although compared to Theorem \ref{as}, there are some extra assumptions (\eqref{T>} and \eqref{assu2} and the last term in \eqref{assu1}) for Theorem \ref{Lp}, we can check that the example in \cite[Section 2]{LY3} also satisfies \eqref{T>} and \eqref{assu2} by a similar argument used in the proof of the example, see \cite[Section 5]{LY3}; and in that example, $\lambda_2$ can be any positive constant thus the last term in \eqref{assu1} is satisfied.
We stress that the example in \cite{LY3} and the second example in this subsection are of different types. While the former has a continuous, unbounded law of displacement which is independent of the branching law, the latter has a discrete, bounded law of displacement which may depend on the branching law.


The constructions of the two examples in this subsection have the following five steps in common.
\begin{itemize}
\item Consider a two-environment case, that is, there are two elements $\omega,\tilde{\omega}\in\Pi$ such that $\mbfP(\L_1=\omega)=p\in(0,1)$ and $\mbfP(\L_1=\tilde{\omega})=1-p.$
Denote $\mbfP_{\omega}:=\mbfP(\cdot|\L_1=\omega), \mbfP_{\tilde{\omega}}:=\mbfP(\cdot|\L_1=\tilde{\omega})$ and $\mbfE_{\omega}, \mbfE_{\tilde{\omega}}$ the corresponding expectations of $\mbfP_{\omega}, \mbfP_{\tilde{\omega}}$.
In the rest of this subsection we always write $N(\phi), \zeta_i(\phi)$ as $N, \zeta_i$ for simplicity.

\item Suppose that $\mbfP_{\omega}(N\leq 2024)=\mbfP_{\tilde{\omega}}(N\leq 2024)=1.$

\item The distributions of $\omega_1$ and $\omega_2$ can be represented as
$$\forall \theta>0,~\mbfE_{\omega}\left(\sum_{i=1}^{N}e^{-\theta\zeta_i}\right)=\sum_{i\leq l}a_ie^{-\theta b_i}:=L(\theta),~\mbfE_{\tilde{\omega}}\left(\sum_{i=1}^{N}e^{-\theta\zeta_i}\right)=\sum_{i\leq \tilde{l}}\tilde{a}_ie^{-\theta \tilde{b}_i}:=\tilde{L}(\theta),$$
where $b_1<b_2<...<b_l, \tilde{b}_1<\tilde{b}_2<...<\tilde{b}_{\tilde{l}}, a_i>0, \tilde{a}_i>0, l\geq 2, \tilde{l}\geq 2.$ (For example, if $$\mbfP_{\omega}\left(\zeta_1=-1,\zeta_2=\zeta_3=0,\zeta_4=\zeta_5=\zeta_6=2\right)=0.1,$$
$$\mbfP_{\omega}\left(\zeta_1=\zeta_2=-1,\zeta_3=2,\zeta_4=3\right)=0.7,~~\mbfP_{\omega}\left(N=0\right)=0.2,$$then
$$L(\theta)=1.5e^{\theta}+0.2+e^{-2\theta}+0.7e^{-3\theta}.)$$

\item  Suppose that $p\log(\sum^{l}_{i=1}a_i)+(1-p)\log(\sum^{\tilde{l}}_{i=1}\tilde{a}_i)>0$ and $p\log a_1+(1-p)\log\tilde{a}_1<0.$

\item  Suppose that $\min\{\sum^{l}_{i=1}a_i, \sum^{\tilde{l}}_{i=1}\tilde{a}_i\}\geq 1.$

\end{itemize}
After these five common steps, we can give the following two examples with different values of $\max\{a_1,\tilde{a}_1\}.$
\begin{example}\label{e2}
(1) If $\max\{a_1 , \tilde{a}_1\}\geq 1$, then for any $p\in(0,1)$, the BRWre satisfies \eqref{qqq} and  all assumptions in Theorem \ref{Lp} except \eqref{T<}.

(2) If $\max\{a_1 , \tilde{a}_1\}<1$, then we can find $0\leq c_-<c_+\leq 1$ such that the BRWre satisfies all the assumptions in Theorem \ref{Lp} when $p\in(c_-,c_+).$
\end{example}

{\it Proof} First, \eqref{c2'a}-\eqref{c3a} hold with any finite positive constants
$\lambda_1, \lambda_{2}, \lambda_{3}$ and $\lambda_{4}$ because the displacement and branching of $\omega$ and $\tilde{\omega}$ are both bounded.

Second, we see \eqref{sect-2} holds from the following observation. Denote $$\Lambda(\theta):=\log L(\theta)-\theta \frac{L'(\theta)}{L(\theta)},~~~\tilde{\Lambda}(\theta):=\log\tilde{L}(\theta)-\theta\frac{\tilde{L}'(\theta)}{\tilde{L}(\theta)}.$$ By a standard argument \footnote{Note that $\Lambda'(\theta)=-\theta(\log L(\theta))''$ and $\log L(\theta)$ is strictly convex because of $l\geq2.$} we see \begin{eqnarray}\label{exx1}\forall \theta>0~~\Lambda'(\theta)<0,~~ \tilde{\Lambda}'(\theta)<0,\end{eqnarray} and \begin{eqnarray}\label{exx2}\lim_{\theta\ra +\infty}\Lambda(\theta)=\log a_1,~~ \lim_{\theta\ra +\infty}\tilde{\Lambda}(\theta)=\log \tilde{a}_1.\end{eqnarray}
Note that $$\kappa(\theta)-\theta\kappa'(\theta)=p\Lambda(\theta)+(1-p)\tilde{\Lambda}(\theta), \forall \theta\geq 0.$$ Therefore, $p\log(\sum^{l}_{i=1}a_i)+(1-p)\log(\sum^{\tilde{l}}_{i=1}\tilde{a}_i)>0$ and $p\log a_1+(1-p)\log\tilde{a}_1<0$ imply $\kappa(0)>0$ and $\lim_{\theta\ra +\infty}(\kappa(\theta)-\theta\kappa'(\theta))<0$ respectively. Moreover, \eqref{exx1} means $(\kappa(\theta)-\theta\kappa'(\theta))'<0$ for all $\theta>0.$ Then we see there exists $\vartheta>0$ such that $\kappa(\vartheta)=\vartheta\kappa'(\vartheta).$

Third, we verify that $\min\{\sum^{l}_{i=1}a_i, \sum^{\tilde{l}}_{i=1}\tilde{a}_i\}\geq 1$ implies \eqref{T>}.
We remind that conditionally on $\omega$ (resp.$\tilde{\omega}$), $\kappa_1(\theta)=\log L(\theta)$ (resp. $\kappa_1(\theta)=\log \tilde{L}(\theta)$).
Note that
$$\mbfE_{\omega}\left(\sum_{i=1}^{N}\1_{\{\theta \zeta_i+ \kappa_1(\theta)\geq 0\}}\right)\geq \mbfE_{\omega}\left(\sum_{i=1}^{N}\1_{\{\zeta_i=b_l,\theta b_l+ \log L(\theta)\geq 0\}}\right)$$
and $\theta b_l+ \log L(\theta)>\Lambda(0), \forall \theta> 0.$ Then $\Lambda(0)=\log(\sum^{l}_{i=1}a_i)\geq 0$ means that
$$\forall \theta>0 ~~(\text{including}~~ \vartheta),~~\mbfE_{\omega}\left(\sum_{i=1}^{N}\1_{\{\zeta_i=b_l,\theta b_l+ \log L(\theta)\geq 0\}}\right)=\mbfE_{\omega}\left(\sum_{i=1}^{N}\1_{\{\zeta_i=b_l\}}\right)=a_l>0.$$
On the other hand, since the branching and displacement are bounded, it is plain to see that $$\mbfE_{\omega}\left(\sum_{i=1}^{N}\1_{\{\vartheta \zeta_i+ \kappa_1(\vartheta)\geq 0\}}\right)=\mbfE_{\omega}\left(\sum_{i=1}^{N}\1_{\{N\leq \lambda,~\vartheta \zeta_i+ \kappa_1(\vartheta)\in [0,\lambda]\}}\right)$$ for $\lambda$ large enough. 
Hence there exists $\lambda>0$ such that $\mbfE_{\omega}\left(\sum_{i=1}^{N}\1_{\{N\leq \lambda,~\vartheta \zeta_i+ \kappa_1(\vartheta)\in [0,\lambda]\}}\right)>0$. By a same argument we also get  $\mbfE_{\tilde{\omega}}\left(\sum_{i=1}^{N}\1_{\{N\leq \lambda',\vartheta \zeta_i+ \kappa_1(\vartheta)\in [0,\lambda']\}}\right)>0$ for some $\lambda'>0$. So far, we have verified \eqref{T>}.

In the rest of the proof, we treat (1) and (2) in different ways. To complete the proof of (1), we only need to explain why \eqref{qqq} follows from $\max\{a_1 , \tilde{a}_1\}\geq 1.$ Without loss of generality, we assume that $a_1\geq 1.$ Note that $\theta b_1+ \log L(\theta)=\log(\sum_{i=1}^la_ie^{\theta(b_1-b_i)})>\log a_1$ for all $\theta>0$ (including $\vartheta$), and hence there exists $\ee:=(\vartheta b_1+ \log L(\vartheta)+\log a_1)/2>0$ such that
$$(0\leq)\mbfE_{\omega}\left(\sum_{i=1}^{N}\1_{\{\vartheta \zeta_i+\kappa_1(\vartheta) \leq \ee\}}\right)\leq \mbfE_{\omega}\left(\sum_{i=1}^{N}\1_{\{\vartheta b_1+ \log L(\vartheta)\leq \ee\}}\right)=0,$$ which means
$\mbfP\left(\mbfE_{\L}\left(\sum_{i=1}^{N}\1_{\{\vartheta \zeta_i+ \kappa_1(\vartheta)\leq \ee\}} \right)=0\right)\geq p>0,$ i.e., \eqref{qqq} holds.

At last, we verify \eqref{T<} by the assumption $\max\{a_1, \tilde{a}_1\}<1$  in (2). Denote $$\tau(\theta):=\theta b_1+\log L(\theta),~~~\tilde{\tau}(\theta):=\theta \tilde{b}_{1}+\log \tilde{L}(\theta).$$ By a standard argument we see
\begin{eqnarray}\label{exx3+}\lim_{\theta\ra+\infty}\tau(\theta)=\log a_1,~~~ \lim_{\theta\ra+\infty}\tilde{\tau}(\theta)=\log \tilde{a}_1.\end{eqnarray}
Since $\min\{l,\tilde{l}\}\geq 2$, we have $$-b_l< \frac{L'(\theta)}{L(\theta)}<-b_1,~~-\tilde{b}_{\tilde{l}}<\frac{\tilde{L}'(\theta)}{\tilde{L}(\theta)}<-\tilde{b}_{1},~~ \forall \theta\geq 0,$$
and hence
\begin{eqnarray}\label{exx3}\max\{\tau'(\theta),\tilde{\tau}'(\theta)\}<0, ~~\forall \theta\geq 0.\end{eqnarray}
Denote
$$\theta_0:=\begin{cases}\text{the zero point of} ~\tau(\cdot),~ \text{if}~ \tau(0)>0\\0, \text{if}~ \tau(0)\leq 0\end{cases}~\text{and}~ \tilde{\theta}_0:=\begin{cases}\text{the zero point of} ~\tilde{\tau}(\cdot),~ \text{if}~ \tilde{\tau}(0)>0\\0, \text{if}~ \tilde{\tau}(0)\leq 0\end{cases}.$$
Obviously, $\theta_0$ and $\tilde{\theta}_0$ are well-defined because of \eqref{exx3}. Moreover, $\max\{\tau(0),\tilde{\tau}(0)\}>0$ follows from $\kappa(0)>0$ and $\max\{\lim_{\theta\ra+\infty}\tau(\theta), \lim_{\theta\ra+\infty}\tilde{\tau}(\theta)\}<0$ from $\max\{a_1, \tilde{a}_1\}<1$ and \eqref{exx3+}. From the analysis above, we see $\bar{\theta}:=\max\{\theta_0,\tilde{\theta}_0\}\in(0,+\infty).$ Without loss of generality, we assume $\bar{\theta}=\theta_0.$ Note that $\Lambda(\bar{\theta})>\tau(\bar{\theta})=0$ since  $L'(\bar{\theta})/L(\bar{\theta})<-b_1.$ Now we choose $$c_+:=1~~~  {\rm and} ~~~c_{-}:=-\frac{\min\{0,\tilde{\Lambda}(\bar{\theta})\}}{\Lambda(\bar{\theta})-\tilde{\Lambda}(\bar{\theta})},$$ which ensures that $\kappa(\bar\theta)-\bar{\theta}\kappa'(\bar{\theta})=p\Lambda(\bar{\theta})+(1-p)\tilde{\Lambda}(\bar{\theta})>0$ as long as $p\in(c_-,c_+).$
Since the functions $\tau(\cdot)$,~ $\tilde{\tau}(\cdot)$,~ $\Lambda(\cdot)$ and $\tilde{\Lambda}(\cdot)$ all strictly deceases at $[0,+\infty),$ we derive that $\vartheta>\bar{\theta}$ and $\tilde{\tau}(\bar{\theta})\leq \tilde{\tau}(\tilde{\theta}_0)=0$. Therefore, it is true that $\max\{\tau(\vartheta),\tilde{\tau}(\vartheta)\}<\max\{\tau(\bar{\theta}),\tilde{\tau}(\bar{\theta})\}=0.$

Recall that $\mbfP_{\omega}(N\leq 2024)=1$ and $\tau(\theta):=\theta b_1+ \log L(\theta).$ Then $\tau(\vartheta)<0$ implies that there exists $\lambda_{51}\leq -1$ such that
\begin{eqnarray}\label{bu1}
\mbfE_{\omega}\left(\sum_{i=1}^{N}\1_{\{N\leq |\lambda_{51}|, \vartheta \zeta_i+ \kappa_1(\vartheta)\leq [\lambda_{51},\lambda^{-1}_{51}]\}}\right)&=&\mbfE_{\omega}\left(\sum_{i=1}^{N}\1_{\{\vartheta \zeta_i+ \kappa_1(\vartheta)\leq \lambda^{-1}_{51}\}}\right)\no
\\&\geq& \mbfE_{\omega}\left(\sum_{i=1}^{N}\1_{\{\zeta_i=b_1,\vartheta b_1+ \log L(\vartheta)\leq \lambda^{-1}_{51}\}}\right)\no
\\&=& \mbfE_{\omega}\left(\sum_{i=1}^{N}\1_{\{\zeta_i=b_1\}}\right)\no
\\&=& a_1(>0).\no
\end{eqnarray}
By a same argument we also get  $$\exists \lambda_{52}\leq -1,~~ \mbfE_{\tilde{\omega}}\left(\sum_{i=1}^{N}\1_{\{N\leq |\lambda_{52}|, \vartheta \zeta_i+ \kappa_1(\vartheta)\leq [\lambda_{52},\lambda^{-1}_{52}]\}}\right)>0,$$ thus \eqref{T<} holds by taking $\lambda_5:=\min(\lambda_{51},\lambda_{52})$.\qed

We mention that the assumption $\min\{\sum^{l}_{i=1}a_i, \sum^{\tilde{l}}_{i=1}\tilde{a}_i\}\geq 1,$ which is set to ensure \eqref{T>}, is not a necessary condition for our example. Here we impose this assumption for the sake of a simpler description in Example \ref{e2}. In fact, for the case of $\min\{\sum^{l}_{i=1}a_i, \sum^{\tilde{l}}_{i=1}\tilde{a}_i\}<1,$ while we can still show that \eqref{T>} holds by a discussion similar to that of verifying $\eqref{T<}$ in the previous paragraph, the expressions of the range of $p$ (i.e.,~$p\in(0,1)$ in (1) and $p\in(c_-,c_+)$ in (2)) will be more complex---the range of $p$ will depend on the sequences $\{a_i\}_{i\leq l}$, $\{\tilde{a}_i\}_{i\leq \tilde{l}},$ $\{b_i\}_{i\leq l},$ $\{\tilde{b}_i\}_{i\leq\tilde{l}}$ and will not be empty as long as we set a proper relationship among these four sequences.

At last, we stress again that the example in \cite[Section 2]{LY3}, whose construction differs totally from that of Example \ref{e2}(2), also satisfies all the assumptions in Theorem \ref{Lp}. Of course, more examples of different types (which satisfy all the assumptions in Theorem \ref{Lp}) could be constructed according to the readers' interests.

\section{Preliminary}

In this section, we give two important tools as the preliminary to prove Theorem \ref{Lp}. One is the many-to-one formula---a kind of measure transformation which has been widely applied in the research of branching random walk. The other is a corollary of strong approximation, which help us to estimate certain kinds of trajectories' distributions precisely.

\subsection{Many-to-one formula}

The many-to-one formula can be traced down to the early works of Peyrie\'{e}re \cite{P1974}~and Kahane and Peyrie\'{e}re \cite{KP1976}. Many variations of this result have been introduced, see e.g. [8]. In this article we need a time-inhomogeneous and bivariate version of many-to-one formula, which has been introduced in \cite{LY3}. For the sake of readability, we retell it in this subsection.

Let $\tau_{n,\L}$ be a random probability measure on $\bfR\times\bfN$ such that for any $x\in\bfR, A\in\bfN,$ we have
\begin{align}\label{m1}\tau_{n,\L}((-\infty,x]\times[0,A])=\frac{\mbfE_{\L}\big(1_{\{N(u)\leq A\}}\sum^{N(u)}_{i=1}1_{\{\zeta_i(u)\leq x\}}e^{-\vartheta \zeta_i(u)}\big)}{\mbfE_{\L}\big(\sum^{N(u)}_{i=1}e^{-\vartheta \zeta_i(u)}\big)},~~|u|=n-1,\end{align}
where $\vartheta$ has been introduced in \eqref{sect-2}.
Hence we can see that the randomness of $\tau_{n,\L}$ comes entirely from $\L_n.$ Moreover, since $N(u)$ only takes values on $\bfN,$ we have
\begin{align*}\tau_{n,\L}(\bfR\times([0,+\infty)\setminus \bfN))=0, ~~~\rm{\mathbf{P}-a.s.}\end{align*}
Under the quenched law $\mbfP_{\L},$ we introduce a series of independent two-dimensional random vectors $\{X_n, \xi_n\}_{n\in\bfN^+}$ whose distributions are $\{\tau_{n,\L}\}_{n\in\bfN^+}.$

Define \beqlb\label{shift0}\chi_n:=\sum_{i=1}^{n}X_i,~ ~T_0:=0,~ ~T_n:=K_n+\vartheta \chi_n, ~\forall n\in\bfN^+.\eeqlb
The $\{T_n\}$ is usually called the associated walk (to the BRWre we consider). We see that under $\mbfP_{\L}$, both $\chi_n$ and $T_n$ are the sums of $n$ independent random variables, and the laws of $X_j$ and $T_j-T_{j-1}$ only depend on $\L_j.$ Moreover, for any measurable function $f,$ 
$\{\mbfE_{\L}f(T_j-T_{j-1})\}_j$ and $\{\mbfE_{\L}f(X_j)\}_j$) are sequences of i.i.d. random variables under $\mbfP$ as the random environment $\L$ is i.i.d. 
\footnote{In fact, $\{T_n\}$ is the random walk with random environment in time studied in \cite{Lv201802}. This model is also a topic discussed in the present paper, see Section 4.} The many-to-one formula shows the connection between BRWre and the associated walk.

\begin{lemma}\label{mto}
(Many-to-one \cite[Lemma 3.1]{LY3}) For any~$n\in\bfN^+$, a positive sequence~$\{A_i\}_{i\in \bfN^+}$ and a measurable function~$f:\bfR^n\rightarrow [0, +\infty),$ we have
\begin{eqnarray}\label{mto0}
&&\mbfE_{\L}\left[\sum_{|u|=n}f(V(u_i),1\leq i\leq n)\1_{\{N(u_{i-1})\leq A_i,1\leq i\leq n\}}\right]\nonumber
\\&&~~~~~~~~~~~~~~~~=\mbfE_{\L}\left[e^{T_n}f(\chi_i,1\leq i\leq n)\1_{\{\xi_{i}\leq A_i,1\leq i\leq n\}}\right],~~~{\rm\mathbf{P}-a.s.}\end{eqnarray}\end{lemma}

This lemma is consistent with \cite[Lemma 2.2]{M2015a} when $A_i=+\infty, \forall i\in\bfN.$ We mention that the distribution of $\{\chi_n\}$ under $\mbfP_{\L}$ is the same as the distribution of the ``spine" in the spinal decomposition theorem \cite[Proposition 2.1]{M2015a} under a probability derived by a size-biased construction, see \cite[Section 2]{M2015a} for details.

Now we are ready to prove Proposition \ref{rem1+}.

\noindent{\it Proof of Proposition \ref{rem1+}}
According to \eqref{mto0}, we have
\begin{align}\label{prop1}\mbfE_{\L}T_1=\kappa_1(\vartheta)-\vartheta\kappa'_1(\vartheta),~~  e^{\kappa_1(\vartheta-\lambda)-\kappa_1(\vartheta)}=\mbfE_{\L}(e^{\lambda \chi_1}),~~e^{\kappa_1(\vartheta+\lambda)-\kappa_1(\vartheta)}=\mbfE_{\L}(e^{-\lambda \chi_1}).\end{align}
(It is worthwhile to remind that the truth of \eqref{prop1} and the forthcoming \eqref{bu2} and \eqref{bu3} have nothing to do with whether the random environment is degenerate or not.)

We stress that in this proof, for any measurable $X$, $\mbfE_{\L}X$ has degenerate distribution as the random environment is degenerate. Hence for any $\theta\geq 0$, $\kappa_1(\theta)$ can be seen as a constant. 
Note that $\kappa_1(\theta)$ is a convex function. \eqref{sect-2} means that
$\kappa_1(0)>0$ and there exist $\vartheta, \lambda>0$ such that $\kappa_1(\vartheta)=\vartheta\kappa'_1(\vartheta)$ and $|\kappa_1(\theta)|<+\infty$ for any $\theta\in[\vartheta-\lambda, \vartheta+\lambda]$. 
Therefore, combining with \eqref{prop1}, we see
\begin{align}\label{prop2}\mbfE_{\L}T_1=0,~~~ \mbfE_{\L}(e^{\frac{\lambda}{\vartheta}|T_1|})\leq \mbfE_{\L}(e^{|\lambda \chi_1+\frac{\lambda}{\vartheta}\kappa_1(\vartheta)|})\leq e^{\frac{\lambda}{\vartheta}|\kappa_1(\vartheta)|}\mbfE_{\L}(e^{\lambda |\chi_1|})<+\infty.\end{align}
Moreover, \cite[Proposition 2.1]{LY3} tells that \eqref{sect-2} implies $\kappa''(\vartheta)>0.$  According to \eqref{mto0}, we have
\begin{eqnarray}\label{bu2}\vartheta^2\kappa_1''(\vartheta)=\mbfE_{\L}[(T_1-\mbfE_\L T_1)^2].\end{eqnarray}
Note that in this proof, it is true that $\kappa''(\vartheta)=\kappa''_1(\vartheta)=\vartheta^{-2}\mbfE_{\L}(T_1^2)$ and thus $\mbfE_{\L}(T_1^2)>0.$ Combining $\mbfE_{\L}(T_1^2)>0$ with \eqref{prop2}, one can find $\lambda_5\leq -1$ such that \begin{align}\label{prop3}\mbfP_\L(T_1\in[\lambda_5,\lambda_5^{-1}])\mbfP_\L(T_1\in[|\lambda_5|^{-1}, |\lambda_5|])>0.\end{align}
According to \eqref{mto0}, we have
\begin{eqnarray}\label{bu3}\mbfP_{\L}(\xi_1>x)=\frac{\mbfE_{\L}\Bigg(\1_{\{N(\phi)>x\}}\sum\limits_{j=1}^{N(\phi)}e^{-\vartheta\zeta_{j}(u)}\Bigg)}{\mbfE_{\L}\Bigg(\sum\limits_{j=1}^{N(\phi)}e^{-\vartheta\zeta_{j}(\phi)}\Bigg)}.\end{eqnarray}
Therefore,
\begin{align}\label{prop3}\mbfP_{\L}(\xi_1>x)\ra 0~~~ \text{as}~~~ x\ra +\infty\end{align} follows from $\mbfE_{\L}\Bigg(\sum\limits_{j=1}^{N(\phi)}e^{-\vartheta\zeta_{j}(\phi)}\Bigg)=e^{\kappa_1(\vartheta)}<+\infty.$

On the other hand, since $\1_{N(\phi)\leq |\lambda_5|}\sum_{i=1}^{N(\phi)}\1_{\{\cdot\}}\leq |\lambda_5|,$ we see \eqref{T>} is equivalent to
$$\mathbf{E}\left(\left|\log\mbfE_{\L}\left(\1_{N(\phi)\leq |\lambda_5|}\sum_{i=1}^{N(\phi)}\1_{\{\vartheta \zeta_i(\phi)+ \kappa_1(\vartheta)\in [0, |\lambda_5|]\}}\right)\right|^{\lambda_6}\right)<+\infty.$$
Note that \eqref{mto0} implies that
$$\mbfE_{\L}\left(\1_{N(\phi)\leq |\lambda_5|}\sum_{i=1}^{N(\phi)}e^{-\vartheta \zeta_i(\phi)-\kappa_{1}(\vartheta)}\1_{\{\vartheta \zeta_i(\phi)+ \kappa_1(\vartheta)\in [0, |\lambda_5|]\}}\right)=\mbfP_{\L}(T_1\in[0,|\lambda_5|], \xi_1\leq |\lambda_5|),$$
hence assumption \eqref{T>} is equivalent to 
\begin{eqnarray}\label{prop1+} \mbfE(|\log\mbfP_{\L}(T_1\in[0,|\lambda_5|], \xi_1\leq |\lambda_5|)|^{\lambda_6})<+\infty.\end{eqnarray}
For the same reason, assumption \eqref{T<} is equivalent to
\begin{eqnarray}\label{prop2+}~~ \mbfE(|\log\mbfP_{\L}(T_1\in[\lambda_5,\lambda^{-1}_5], \xi_1\leq |\lambda_5|)|^{\lambda_6})<+\infty.\end{eqnarray}
(We remind that the statement ``\eqref{T>} and \eqref{T<} are equivalent to \eqref{prop1+} and \eqref{prop2+} respectively" is true even though the random environment is not degenerate.)
Therefore, in the context of degenerate environment, the statement ``both \eqref{T<} and \eqref{T>} hold" are equivalent to
\begin{eqnarray}\label{prop5}\exists \lambda_5\leq-1,~~ \mbfP_{\L}(T_1\in[0,|\lambda_5|], \xi_1\leq |\lambda_5|)\mbfP_{\L}(T_1\in[\lambda_5,\lambda^{-1}_5], \xi_1\leq |\lambda_5|)>0.\end{eqnarray}
If $\mbfP_{\L}(T_1\in[0,|\lambda_5|], \xi_1\leq |\lambda_5|)=0$ for any $\lambda_5\leq-1$, which causes that \eqref{prop5} fails to hold, then
$$\forall \lambda_5\leq-1,~\mbfP_{\L}(T_1\notin[0,|\lambda_5|])+\mbfP_{\L}(\xi_1> |\lambda_5|)\geq \mbfP_{\L}(\{T_1\notin[0,|\lambda_5|]\}\cup \{\xi_1> |\lambda_5|\})=1.$$
Combining with \eqref{prop3} and the second term in \eqref{prop2}, the display above means $\mbfP_{\L}(T_1<0)=1,$ which is a contradiction to $\mbfE_{\L}T_1=0$. If $\forall \lambda_5\leq-1,~\mbfP_{\L}(T_1\in[\lambda_5,\lambda^{-1}_5], \xi_1\leq |\lambda_5|)=0$ causes that \eqref{prop5} fails to hold, then we can get $\mbfP_{\L}(T_1\geq \lambda^{-1}_5)=1$ for any $\lambda_5<-1,$ which contradicts the facts that $\mbfE_{\L}(T_1^2)>0$ and $\mbfE_{\L}T_1=0$. As a result, we see \eqref{prop5} holds and thus both \eqref{T<} and \eqref{T>} hold.\qed

\subsection{Strong approximation}

The other important tools used frequently in the forthcoming proof (in the next section) are the celebrated Sakhanenko's strong approximation theorem and its corollary, which will be stated in the following Theorem \uppercase\expandafter{\romannumeral1} and Corollary \ref{sat}.
Let~$V_1,V_2,\ldots,V_n,\ldots$ be a sequence of independent random variables satisfying $\forall j, \bfE(V_j)=0$ and $\bfE(V^2_j)<+\infty$. Denote $\mathcal{D}_k:=\sum_{i=1}^k\bfE(V^{2}_i).$ Introduce a random broken line $\mathcal{V}(s), s\in\bfR^+$ such that $\mathcal{V}(0)=0,$ $\mathcal{V}(\mathcal{D}_k)=\sum_{i=1}^{k}V_i, k\in\bfN^+$ and $\mathcal{V}(\cdot)$ is linear, continuous on each interval $[\mathcal{D}_{k-1},\mathcal{D}_{k}].$ 
The following theorem is known as the Sakhanenko's strong approximation theorem with power moment.

\noindent\emph{{\bf Theorem \uppercase\expandafter{\romannumeral1}~(Sakhanenko, \cite[Theorem 1]{Sak2006})}
For any $\beta\geq 2,$ there exists a standard Brownian motion $B$ such that
\begin{eqnarray}\label{sak06}\forall x>0,~~~ \bfP\left(\sup_{s\leq \mathcal{D}_n}\big|\mathcal{V}(s)-B_s\big|\geq 2C_0\beta x\right)\leq
\frac{\sum_{k=1}^n\bfE(|V_k|^{\beta})}{x^{\beta}},\end{eqnarray}}
where $C_{0}$ is an absolute constant.

\noindent\emph{{\bf Theorem \uppercase\expandafter{\romannumeral2}~(Cs\"{o}rg\H{o}~and~R\'{e}v\'{e}sz, \cite[Lemma 1]{CR1979})}} For a standard Brownian motion $B$ and a constant $D_1>2$, there exists a constant $D_2\in(0,+\infty)$ (depending only on $D_1$) such that
$$\forall x>0,~~~t>0,~~~ \bfP\left(\sup_{0\leq s\leq t}|B_s|\geq x\right)\leq D_2e^{-\frac{x^2}{D_1t}}.$$
\begin{defn}
Let $V$ be a random variable with mean $0.$ We call $(m,l)$ a space-time adapted (STA) couple at level $(\beta, \iota)$ with respect to $V$ if
$m>0, l\in\bfN^+, \beta\geq 2, \iota>2,\mbfE(|V|^{\beta})<+\infty$ and $\iota\mbfE(|V|^{2})lm^{-2}\log(lm^{-\beta})\geq-1.$
\end{defn}
\begin{cor}\label{sat}
Let $\{V_i\}$ be a sequence of i.i.d. copies of $V$. If $(m,l)$ is a STA couple at level $(\beta, \iota)$ w.r.t. $V_1$, then we can find a constant $C$ depending only on $\beta, \iota$ and $\mbfE(|V|^{\beta})$ such that
$$\bfP\left(\max_{i\leq l}\left|\sum_{k=1}^iV_k\right|\geq m\right)\leq C\frac{l}{m^\beta}.$$
\end{cor}

(Note that the corollary above can not be obtained by Doob's inequality. Since $\{\sum_{k=1}^iV_k\}_i$ is a martingale, Doob's inequality tells that $\bfP\left(\max_{i\leq l}\left|\sum_{k=1}^iV_k\right|\geq m\right)\leq \frac{\mbfE(|\sum_{k=1}^lV_k|^{\beta})}{m^\beta}.$ But $\lim_{l\ra+\infty}l^{-1}\mbfE(|\sum_{k=1}^lV_k|^{\beta})=+\infty$ as long as $\beta>2.$)

{\it Proof of Corollary \ref{sat}} Note that there exists $\iota_0\in(0,1)$ such that $2+\iota_0=(1-\iota_0)\iota$ because of $\iota>2$. Then $D^{-1}_1t^2+\mbfE(|V|^{2})lm^{-2}\log(lm^{-\beta})\geq 0$ holds when we take $t=\sqrt{1-\iota_0}$ and $D_1=2+\iota_0.$ 
From Theorem \uppercase\expandafter{\romannumeral1} we can find a  Brownian motion $B$ such that
$$\mathcal{P}_1:=\bfP\left(\max_{i\leq l}\left|\left(\sum_{k=1}^iV_k\right)-B_{i\bfE(V^2_1)}\right|\geq (1-t)m\right)\leq \left(\frac{2C_0\beta}{1-t}\right)^{\beta}\bfE(|V_1|^{\beta})\frac{l}{m^\beta}.$$
Theorem \uppercase\expandafter{\romannumeral2} tells that
$$\mathcal{P}_2:=\bfP\left(\max_{i\leq l}\left|B_{i\bfE(V^2_1)}\right|\geq tm\right)\leq D_2\exp\left\{-\frac{t^2m^2}{D_1\bfE(V^2_1)l}\right\}$$
and hence $\mathcal{P}_2\leq D_2\frac{l}{m^{\beta}}.$
Note that
\begin{eqnarray}\label{citeprop}\bfP\left(\max_{i\leq l_n}\left|\sum_{k=1}^iV_k\right|\geq m_n\right)\leq \mathcal{P}_1+\mathcal{P}_2,\end{eqnarray} and recall that $D_2$ and $t$ are only determined by $\iota$, hence we complete the proof by taking $C:=\left(\frac{2C_0\beta}{1-t}\right)^{\beta}\bfE(|V_1|^{\beta})+D_2.$


\section{Proof of Theorem \ref{Lp}}
We divide the proof into two parts. In the first subsection, we state the idea and lead readers to see what is the most challenging part in the proof, then we will prove the most challenging part in the second subsection.

\subsection{A lower bound of survival probability}
Let us first give some classical conclusions as lemmas, which inspire us to develop the approach to Theorem \ref{Lp}.
\begin{lem}\label{l1}
For a non-negative random sequence $\{g_n\}$ and constants $c\geq 0, p\geq 1,$ if we have $\bfE|g^q_n-c^q|\ra 0, \forall q\in(0,p],$ then  $\bfE(|g_n-c|^p)\ra 0$ \footnote{We usually write $\bfP, \bfE$ for the probability and the corresponding expectation when there is no random environment involved.}.
\end{lem}
{\bf Proof:} It follows the facts that $|g_n-c|^p\leq \max(2^{p-2},1)|g_n-c|(g^{p-1}_n+c^{p-1})$ and
$$(g_n-c)(g^{p-1}_n+c^{p-1})=g^p_n-c^p-c(g^{p-1}_n-c^{p-1})+c^{p-1}(g_n-c).$$
~~~~~~~~~~~~~~~~~~~~~~~~~~~~~~~~~~~~~~~~~~~~~~~~~~~~~~~~~~~~~~~~~~~~~~~~~~~~~~~~~~~~~~~~~~~~~~~~~~~~~~~~~~~~~~~~~~~~~\qed

The next lemma is known as the Vitali convergence theorem, which has been stated in many textbooks in various forms. Here we state a refined version of \cite[Exercise 7.17]{R2013}.
\begin{lem}\label{l2}
(Vitali convergence theorem) If $X_n$ converges to a constant $c$ in probability $\bfP$ and $\{X_n\}$ is uniformly integrable, then $\bfE|X_n-c|\ra 0.$
\end{lem}
{\bf Proof:} Recall that the uniformly integrability means that for any $\epsilon>0,$ we can find a constant $M>0$ such that $\sup_n\bfE(|X_n|\1_{|X_n|\geq M})<\epsilon.$ Therefore, this lemma immediately follows from the next two inequalities:
$$\bfE|X_n-c|\leq \bfE(|X_n|\1_{|X_n|\geq M})+|c|\bfP(|X_n|\geq M)+\bfE(|X_n-c|\1_{|X_n|< M})$$
and
$$\forall \epsilon>0, ~~\bfE(|X_n-c|\1_{|X_n|< M})\leq \epsilon\bfP\left(|X_n-c|\leq \epsilon\right)+(M+|c|)\bfP\left(|X_n-c|> \epsilon\right).$$
\begin{lem}\label{l3}
(\cite[Exercise 7.19]{R2013}) If there exists $\epsilon>0$ such that $\sup_n\bfE\left(|X_n|^{1+\epsilon}\right)<+\infty,$ then $\{X_n\}$ is uniformly integrable.
\end{lem}


Now we start the first half of the proof of Theorem \ref{Lp}: from the conclusion in Theorem \ref{Lp} to \eqref{goal}.

{\bf Proof of Theorem \ref{Lp}: first half} ~~First we see the
 assumptions in Theorem \ref{inp} are totally contained in the assumptions in Theorem \ref{Lp}, which means that \eqref{b1} and \eqref{b2} hold in probability. Denote
$$A_n:=-\frac{\log\mbfP_{\L}(Y_n>0)}{\sqrt[3]{n}}.$$ Note that $A_n$ converges to a constant in probability means that for any $p\geq 1,$ $A^p_n$ also converges to a constant in probability. Therefore, if $\{A^p_n\}$ is uniformly integrable (and thus $\{A^q_n\}$ is uniformly integrable for $q\in(0,p]$), then according to Lemma \ref{l2}, we see for any $q\in(0,p],$ $\bfE|A_n^q-b^q_1|\ra 0$ (resp. $\bfE|A_n^q-b^q_2|\ra 0$) when $\alpha=\frac{1}{3}, d\in(0,d_c)$ (resp. $\alpha\in(0,\frac{1}{3}), d>0$). According to Lemma \ref{l1}, we see $A_n\ra b_1, L^p$ (resp. $A_n\ra b_2, L^p$) if $\bfE|A_n^p-b^p_1|\ra 0$ (resp. $\bfE|A_n^p-b^p_2|\ra 0$).  Therefore, Lemma \ref{l3} tells that if we can show \begin{eqnarray}\label{An}\exists \varepsilon>0, ~~~\varlimsup_{n\ra +\infty}\mbfE(A_n^{p+\varepsilon})<+\infty,\end{eqnarray}
then we complete the proof of Theorem \ref{Lp}.
Therefore, what we should do next is to find the lower bound of $\mbfP_{\L}(Y_n>0).$
Note that 
$$Y_n\geq \hat{Y}_{n}:=\sharp\left\{|u|=n:\forall i\leq n, V(u_i)\in[-n^{1/3}-\vartheta^{-1} K_i, \varphi_{\L}(i)],~N(u_{i-1})\leq e^{n^{1/3}}\right\},$$
hence \begin{eqnarray}\label{csin}\mbfP_{\L}(Y_n>0)\geq \mbfP_{\L}(\hat{Y}_{n}>0)=\mbfP_{\L}(\hat{Y}_{n}\geq 1)\geq \frac{[\mbfE_{\L}(\hat{Y}_{n}
)]^2}{\mbfE_{\L}(\hat{Y}^2_{n})}, ~~{\rm \mathbf{P}-a.s.,}\end{eqnarray}
where the last inequality is because of the H\"{o}lder's inequality.
From the same argument (which is conventionally called the second moment method) used in \cite[(6.10)-(6.16)]{LY3}, 
we derive that 
\begin{eqnarray}\label{smm}\mbfE_{\L}(\hat{Y}^2_{n})\leq\mbfE_{\L}(\hat{Y}_{n})\left(1+(e^{n^{1/3}}-1)\sum\limits_{j=0}^{n-1}\sup_{\substack{|v|=j \\V(v)\in\mbr}}\mbfE_{\L}\Big[Z_n^{v}(\Theta)\Big]\right), ~~{\rm \mathbf{P}-a.s.,}\end{eqnarray}
where $Z^v_n(\Theta):=\1_{\{v\in\Theta\}}\left(\sum_{|u|=n,u>v}\1_{\{u\in\Theta\}}\right)$ and
$$\Theta:=\left\{u\in \mathbf{T}:|u|\leq n,~~\forall 0\leq i\leq |u|,~~N(u_{i-1})\leq e^{n^{1/3}}, V(u_i)\in[-n^{1/3}-\vartheta^{-1} K_{i},~\varphi_{\L}(i)]\right\}.$$
(In fact, we can obtain \eqref{smm} by just redefining $\Theta$ in \cite{LY3} as the display above and proceeding via the steps (6.10)-(6.16) in \cite{LY3}.)
Applying Lemma \ref{mto} we see
\begin{eqnarray}\label{rough}\sup_{\substack{|v|=j\\V(v)\in \bfR}}\mbfE_{\L}\Big[Z_n^{v}(\Theta)\Big]&=&\sup_{\substack{|v|=j\\V(v)\in \bfR}}\mbfE_{\L}\Bigg[\sum_{\substack{u_j=v,|u|=n}}\1_{\{\forall i\leq n-j,~V(u_{j+i})+\vartheta^{-1} K_{i+j}\in[-n^{1/3},~d(i+j)^{\alpha}]\}}\Bigg]\no
\\&\leq&\sup_{\substack{|v|=j\\V(v)\in \bfR}}\mbfE_{\L}\Bigg[\sum_{\substack{|u|=n}}\1_{\{V(u)+\vartheta^{-1} K_{n}\leq dn^{\alpha}\}}\Bigg]\no
\\&=&\mbfE_{\L}\Bigg[\sum_{\substack{|u|=n}}\1_{\{V(u)+\vartheta^{-1} K_{n}\leq~dn^{\alpha}\}}\Bigg]\no
\\&=&\mbfE_{\L}\Bigg[e^{\vartheta \chi_n+K_n}\1_{\{\chi_n+\vartheta^{-1} K_{n}\leq dn^{\alpha}\}}\Bigg]\no
\\&=&\mbfE_{\L}\Bigg[e^{T_n}\1_{\{T_{n}\leq \vartheta dn^{\alpha}\}}\Bigg]\no
\\&\leq&e^{\vartheta dn^{\alpha}},~~  {\rm \mathbf{P}-a.s.}\end{eqnarray}
(The estimate above is very rough indeed. But we point out that it is enough for our aim because in this step, we only need to prove that $\sup_{|v|=j,V(v)\in\mbr}\mbfE_{\L}\left[Z_n^{v}(\Theta)\right]$ can be bounded by a constant from above.
In fact, a more precise random upper bound of $n^{-\frac{1}{3}}\log\sup_{|v|=j,V(v)\in\mbr}\mbfE_{\L}\left[Z_n^{v}(\Theta)\right]$ had been obtained in the proof of Theorem \ref{as}. Thanks to the conclusion in Theorem \ref{as}, in this proof we only need to give such a rough estimate on $\sup_{|v|=j,V(v)\in\mbr}\mbfE_{\L}\left[Z_n^{v}(\Theta)\right]$.) 

From \eqref{rough} we see whether $\alpha=\frac{1}{3}, d\in(0,d_c)$ or $\alpha\in(0,\frac{1}{3}), d\geq 0,$ it is true that  \begin{eqnarray}\label{rough'}\sup_{|v|=j,V(v)\in\mbr}\mbfE_{\L}\Big[Z_n^{v}(\Theta)\Big]\leq e^{(\vartheta d_c+1)n^{1/3}},~~  {\rm \mathbf{P}-a.s.} \end{eqnarray} for $n$ large enough. Combining with \eqref{csin} and \eqref{smm} we see
\begin{eqnarray}\label{rough'}\mbfP_{\L}(Y_n>0)&\geq&\frac{\mbfE_{\L}(\hat{Y}_{n})}{1+(e^{n^{1/3}}-1)\sum\limits_{j=0}^{n-1}\sup_{\substack{|v|=j \\V(v)\in\mbr}}\mbfE_{\L}\Big[Z_n^{v}(\Theta)\Big]}\no
\\&\geq& \frac{\mbfE_{\L}(\hat{Y}_{n})}{e^{(\vartheta d_c+3)n^{1/3}}},~~~{\rm \mathbf{P}-a.s.}\end{eqnarray}
Using Lemma \ref{mto} once again, we get\footnote{In this proof, we always write $\mbfP_{\L}(\cdot|T_0=0)$ as $\mbfP_{\L}(\cdot)$ for simplicity if no confusion may arise.}
\begin{eqnarray}
\mbfE_{\L}(\hat{Y}_{n})&=&\mbfE_{\L}\left(e^{T_{n}}\1_{\{\forall 1\leq i\leq n, ~T_i\in[-\vartheta n^{1/3},~\vartheta di^{\alpha}],~\xi_{i}\leq e^{n^{1/3}}\}}\right)\nonumber
\\&\geq& e^{-\vartheta n^{1/3}}\mbfP_{\L}\left(\forall_ {1\leq i\leq n},~~\xi_{i}\leq e^{n^{1/3}},~T_i\in[-\vartheta  n^{1/3},~\vartheta di^{\alpha}]\right),\no
\\&\geq& e^{-\vartheta n^{1/3}}\mbfP_{\L}\left(\forall_ {1\leq i\leq n},~~\xi_{i}\leq e^{n^{1/3}},~T_i\in[-\vartheta  n^{1/3},0]\right).\no
\end{eqnarray}
Combining with \eqref{csin} we can see whether $\alpha=\frac{1}{3}, d\in(0,d_c)$ or $\alpha\in(0,\frac{1}{3}), d\geq 0$,
$$A_n\leq \vartheta -n^{-1/3}\log \mbfP_{\L}\left(\forall_{1\leq i\leq n},~\xi_{i}\leq e^{n^{1/3}},~T_i\in[-\vartheta  n^{1/3},0]\right)+(\vartheta d_c+3).$$
Therefore, to show \eqref{An}, we only need to prove that
\begin{eqnarray}\label{Bn}\exists \varepsilon>0, \varlimsup_{n\ra +\infty}\mbfE\left(\left|n^{-1/3}\log \mbfP_{\L}\left(\forall_{1\leq i\leq n},~\xi_{i}\leq e^{n^{1/3}},~T_i\in[-\vartheta n^{1/3},0]\right)\right|^{p+\varepsilon}\right)<+\infty.~\end{eqnarray}
Denote $\lfloor x\rfloor:=\sup\{y\in\bfN, y\leq x\}$ and $\lceil x\rceil:=\inf\{y\in\bfN, y\geq x\}$. For any $R\in(0, \frac{\vartheta}{2|\lambda_5|})$ and $n$ large enough, by Markov property we have
 \begin{eqnarray}\label{mog.2}
 &&\mbfP_\L
\left(\forall_{0\leq i\leq n}~T_{i}\in \left[-\vartheta n^{1/3},0\right],~\xi_{i}\leq e^{n^{1/3}}\right)\nonumber
\\&\geq&\mbfP_\L
(\forall_{0\leq i\leq \lceil Rn^{1/3}\rceil}~T_{i}\in [\lambda_5 i,i/\lambda_5],~\xi_{i}\leq e^{n^{1/3}}|T_{0}=0)\nonumber
\\&\times&\inf_{x\in[2\lambda_5Rn^{1/3},Rn^{1/3}/\lambda_{5}]}\mbfP_\L
\left(\forall_{\lceil Rn^{1/3}\rceil\leq i\leq n}~T_{i}\in [-\vartheta n^{1/3},0],~\xi_{i}\leq e^{n^{1/3}}\Bigg|T_{\lceil Rn^{1/3}\rceil}=x\right)\nonumber
\\&\geq&\prod_{m=0}^{\lceil Rn^{1/3}\rceil-1}\mbfP_\L
(T_{m+1}\in [\lambda_{5},1/\lambda_{5}],\xi_{m+1}\leq |\lambda_5||T_{m}=0)\nonumber
\\&\times&\prod_{j=0}^{\lceil n^{1/3}\rceil}\inf_{x\in[2\lambda_5 Rn^{1/3},Rn^{1/3}/\lambda_{5}]}\mbfP_\L
\left(\begin{aligned}\forall_{R_{j,n}\leq i\leq R_{j+1,n}}~T_{i}\in \left[-\vartheta n^{1/3},0\right],~~~\\~\xi_{i}\leq e^{n^{1/3}},T_{R_{j+1,n}}\in[2\lambda_5 Rn^{1/3},Rn^{1/3}/\lambda_{5}]\end{aligned}\Bigg|T_{R_{j,n}}=x\right)\no
\\&:=&\prod_{m=0}^{\lceil Rn^{1/3}\rceil-1} Z_m\prod_{j=0}^{\lceil n^{1/3}\rceil} Z_{j,n},
\end{eqnarray}
where $R_{j,n}:=\lceil Rn^{1/3}\rceil+\lfloor n^{2/3}\rfloor j$ and $\lambda_5$ is the one in \eqref{T<}.
Note that the random environment is i.i.d., then we have
\begin{eqnarray}\label{Z0n}&&
\mbfE\left(\left|n^{-1/3}\log \mbfP_{\L}\left(\forall_{1\leq i\leq n},~\xi_{i}\leq e^{n^{1/3}},~T_i\in[-\vartheta  n^{1/3},0]\right)\right|^{p+\varepsilon}\right)\no
\\&\leq&\frac{(Rn^{1/3}+n^{1/3}+2)^{p+\varepsilon-1}\left(\lceil Rn^{1/3}\rceil\mbfE\left(|\log Z_0|^{p+\varepsilon}\right)+\lceil n^{1/3}\rceil\mbfE\left(|\log Z_{0,n}|^{p+\varepsilon}\right)\right)}{n^{\frac{p+\varepsilon}{3}}}.
\end{eqnarray}
Recalling \eqref{prop2+} one sees that $\mbfE\left(|\log Z_0|^{p+\varepsilon}\right)<+\infty$ when we choose $\ee\in(0, \lambda_6-p).$
Therefore, \eqref{Bn} will hold as long as
\begin{eqnarray}\label{B+n}\exists \varepsilon>0,~~ \varlimsup_{n\ra +\infty}\mbfE\left(|\log Z_{0,n}|^{p+\varepsilon}\right)<+\infty.\end{eqnarray}
Therefore, from the space-homogeneous property of $\{T_n\}$, \eqref{B+n} will be true if there exists $\ee>0$ such that
\begin{eqnarray}\label{goal}
\forall a,b\in\bfR, ~~b>a>0,~~\varlimsup_{n\ra +\infty}\mbfE\left(\left|\log \hat{p}_{\L}(T;n,1,a,b,a)\right|^{p+\varepsilon}\right)<+\infty,
\end{eqnarray}
where \begin{eqnarray}\label{pt1}\hat{p}_{\L}(T;n,z,a,b,c):=\inf_{|x|\leq a\sqrt{n}}\mbfP_\L
\left(\forall_{i\leq zn, i\in\bfN}~|T_{i}|\leq b\sqrt{n},~|T_{\lfloor zn\rfloor}|\leq c\sqrt{n},~\xi_{i}\leq e^{\sqrt{n}}\big|T_{0}=x\right).~~~~~~~\end{eqnarray}
In conclusion, the proof of Theorem \ref{Lp} will be completed as long as we can prove \eqref{goal}. The proof of \eqref{goal} will be given in the next subsection, which is the heart of the proof of Theorem \ref{Lp}.

 ~~
\qed

\subsection{Proof of \eqref{goal}}
Now we start the second half of the proof of Theorem \ref{Lp}: proof of \eqref{goal}. We divide the proof into three parts.

{\bf Proof of Theorem \ref{Lp}: second half}~~

{\it Part 1.}~Since $\left|\log \hat{p}_{\L}(T;n,1,a,b,a)\right|^{p+\varepsilon}$ is non-negative, it is true that
$$\mbfE\left(\left|\log \hat{p}_{\L}(T;n,1,a,b,a)\right|^{p+\varepsilon}\right)\leq 1+\sum_{m=1}^{+\infty}\mbfP\left(\left|\log \hat{p}_{\L}(T;n,1,a,b,a)\right|^{p+\varepsilon}\geq m\right).$$
Markov property tells that for 
$n$ large enough, we have
\begin{eqnarray}&&\hat{p}_{\L}(T; n,1,a,b,a)\no
\\&\geq& \prod_{i=0}^{n-1}\inf_{|x|\leq a\sqrt{n}}\mbfP_{\L}\left(T_{i+1}\leq a\sqrt{n},\xi_{i+1}\leq |\lambda_5||T_i=x\right)\no
\\&\geq& \prod_{i=0}^{n-1}\min\left(\mbfP_{\L}\left(T_{i+1}\in[0,|\lambda_5|],\xi_{i+1}\leq |\lambda_5||T_i=0\right),\mbfP_{\L}\left(T_{i+1}\in[\lambda_5,0],\xi_{i+1}\leq |\lambda_5||T_i=0\right)\right)\no
\\&\geq& \prod_{i=0}^{n-1}\left(\mbfP_{\L}\left(T_{i+1}\in[0,|\lambda_5|],\xi_{i+1}\leq |\lambda_5||T_i=0\right)\mbfP_{\L}\left(T_{i+1}\in[\lambda_5,0],\xi_{i+1}\leq |\lambda_5||T_i=0\right)\right)\no
\\&:=&\prod_{i=0}^{n-1}(Z_{+,i}Z_{-,i}).\no\end{eqnarray}
Recalling assumption \eqref{assu2} we see
\begin{eqnarray}\label{add}|\log\hat{p}_{\L}(T; n,1,a,b,a)|\leq\left(\left|\sum_{i=0}^{n-1}\log(Z_{+,i}Z_{-,i})\right|^{1/t}\right)^t\leq\left(\sum_{i=0}^{n-1}\left|\log(Z_{+,i}Z_{-,i})\right|^{1/t}\right)^t\end{eqnarray}
because of $t\geq 1.$
According to assumption \eqref{assu2}, we can choose $\varepsilon\in (0, \lambda_6-p)$ such that $\min(\frac{\lambda_0}{2}-\varepsilon,\lambda_1,\frac{\lambda_3}{2})>\frac{\lambda_6}{\lambda_6-p-\varepsilon}.$
Note that \eqref{prop1+} and \eqref{prop2+} mean that 
$$\mbfE(|\log(Z_{+,0}Z_{-,0})|^{\lambda_6})=\mbfE\left(|\log(Z_{+,0}Z_{-,0})|^{\frac{1}{t}\lambda_6t}\right)<+\infty.$$
Denote $\bar{p}:=p+\varepsilon, d_n:=\lceil n^{\frac{\bar{p}}{\lambda_6-\bar{p}}}\rceil$ and hence $d^{1/\bar{p}}_n=O(n^{s})$ for some $s>t$ as $p\geq\lambda_6-\frac{1}{t}.$ 
Moreover,~
note that for $n$ large enough and each $m\geq d_n,$ $(m^{1/(\bar{p}t)},n)$ is a STA couple at level $(\lambda_6t, 3)$ w.r.t. $|\log(Z_{+,0}Z_{-,0})|^{1/t}-z,$ where $z:=\mbfE(|\log(Z_{+,0}Z_{-,0})|^{1/t}).$ We remind that $\{|\log(Z_{+,i}Z_{-,i})|^{1/t}\}_i$ is an i.i.d. sequence. Then
from Corollary \ref{sat} and \eqref{add}, we can find a constant $c_1$ such that
\begin{eqnarray}
\mbfP\left(\left|\log \hat{p}_{\L}(T;n,1,a,b,a)\right|^{p+\varepsilon}\geq m\right)
&\leq&\mbfP\left(\sum_{i=0}^{n-1}(\left|\log(Z_{+,i}Z_{-,i})\right|^{\frac{1}{t}}-z)\geq m^{\frac{1}{\bar{p}t}}-zn\right)\no
\\&\leq&\mbfP\left(\sum_{i=0}^{n-1}(\left|\log(Z_{+,i}Z_{-,i})\right|^{\frac{1}{t}}-z)\geq \frac{1}{2}m^{\frac{1}{\bar{p}t}}\right)\no
\\&\leq& c_1nm^{-\frac{\lambda_6}{\bar{p}}}
\end{eqnarray}
as long as $n$ is large enough.
Hence we have
\begin{eqnarray}
\sum^{+\infty}_{m=d_n}
\mbfP\left(\left|\log \hat{p}_{\L}(T;n,1,a,b,a)\right|^{p+\varepsilon}\geq m\right)
&\leq&c_1n\int_{d_n}^{+\infty}x^{-\frac{\lambda_6}{\bar{p}}}dx\leq 2c_1nd^{1-\frac{\lambda_6}{\bar{p}}}_n<3c_1\no
\end{eqnarray}
for $n$ large enough, which means that
\begin{eqnarray}\label{1<}\varlimsup_{n\ra \infty}\sum^{+\infty}_{m=d_n}
\mbfP\left(\left|\log \hat{p}_{\L}(T;n,1,a,b,a)\right|^{p+\varepsilon}\geq m\right)<+\infty.\end{eqnarray}

{\it Part 2.}~On the other hand, for any event $Q_n,$ we have
$$\mbfP\left(\left|\log \hat{p}_{\L}(T;n,1,a,b,a)\right|^{p+\varepsilon}\geq m\right)\leq \mbfP\left(\left|\log \hat{p}_{\L}(T;n,1,a,b,a)\right|^{p+\varepsilon}\geq m, Q_n\right)+\mbfP\left(Q^c_n\right).$$
Hence
\begin{eqnarray}\label{23<}&&\sum^{d_n}_{m=1}\mbfP\left(\left|\log \hat{p}_{\L}(T;n,1,a,b,a)\right|^{p+\varepsilon}\geq m\right)\no
\\&\leq& d_n\mbfP\left(Q^c_n\right)+\sum^{d_n}_{m=1}\mbfP\left(\left|\log \hat{p}_{\L}(T;n,1,a,b,a)\right|^{p+\varepsilon}\geq m, Q_n\right).\end{eqnarray}
Before introducing the definition of $Q_n$, we give some notation used in the rest of the proof. Denote
$$M_{n}:=\mbfE_\L(T_n), ~~U_{n}:=T_n-\mbfE_\L(T_n),~~\Gamma_n:=\mbfE_\L(U^2_n)=\mbfE_\L(T^{2}_n)-M^2_{n},$$
and$$\psi_n:=\sum_{k=1}^n\mbfE_\L(|U_k-U_{k-1}|^{\lambda_2}).$$
Let us introduce some important constants for the rest of the proof. Thanks to the assumptions $0<a<b, \sigma^2_*>0$ and $\frac{\sigma^2}{\sigma^2_*}< \frac{\lambda_2-2}{\lambda_0-2}$, we can find $\delta, c_3>0$ small enough to satisfy
\begin{eqnarray}\label{range+}c_3\in\left(0,\frac{\sigma^2_*}{2023}\right), \delta\in\left(0,\frac{a}{2023}\right),\max\{0,2a-b\}<a-(1+c_3)\delta,\frac{\sigma^2(1+3c_3)^2}{\sigma^2_*-2c_3}< \frac{\lambda_2-2}{\lambda_0-2},~~~~\end{eqnarray}
which allows us to choose constants $\bar{a}$ and $c_2$ such that
\begin{eqnarray}\label{rangec2-}\bar{a}\in(\max\{0,2a-b\},a-(1+c_3)\delta),~\bar{a}>(1+2c_3)\delta,~\frac{(a-\bar{a}+3c_3\delta)^2}{(b-a-\delta-c_3\delta)^2}\leq \frac{2}{2+c_3\delta},\end{eqnarray}
\begin{eqnarray}\label{rangec2}\frac{c_2(1+2c_3)^2\delta^2}{\sigma^2_*-2c_3}< \lambda_2-2~~\text{and}~~ c_2>\frac{\sigma^2(\lambda_0-2)}{\delta^2}.\end{eqnarray}
Now we give some helpful random sequences as follows.  Define
$$\tau_{_{n,0}}:=0,~~ \tau_{_{n,k}}:=\min\left(\min\{i>\tau_{_{n,k-1}}: |M_{i}-M_{\tau_{_{n,k-1}}}|>\delta\sqrt{n}\}-1, \tau_{_{n,k-1}}+n\right),~~k=1,2,\ldots$$
Denote \begin{eqnarray}\label{rho*}N_n:= \sup\{k\leq n: \tau_{_{n,k}}\leq n\}, ~~~\rho_{_{n,0}}:=0,~~~\rho_{_{n,k}}:=\tau_{_{n,k}}-\tau_{_{n,k-1}}.\end{eqnarray}
The spirit of the construction of $Q_n$ is on which the value of $|\log \hat{p}_{\L}(T;n,1,a,b,a)|$ can not be too large.
 Let $q_n:=\frac{n}{c_2\log n}, ~q^*_n:=\left(\lceil\frac{n}{q_n}\rceil\right)^2$ 
 and define $$Q_n:=H_n\cap J_n\cap \tilde{J}_n \cap I_n,$$ where
$$H_n:=\mbfP\left(\min_{i\leq q^*_n}\rho_{n,i}\geq q_n\right), ~~J_n:=\left\{\max_{i\leq 2n}|\Gamma_{i}-\sigma^2_*i|\leq c_3q_n\right\},~~\tilde{J}_n:=\{\psi_{2n}\leq 3\mbfE(\psi_1)n\},$$
\beqlb\label{tilIn}I_n:=\left\{\max_{|u|\leq n-1}\mbfE_{\L}(N(u)^{1+\lambda_4})\leq e^{\frac{\lambda_4\sqrt{n}}{3}}, \max_{i\leq n}\left[\frac{\vartheta}{\vartheta+\lambda_4}\kappa_i(\vartheta+\lambda_4)-\kappa_i(\vartheta)\right]\leq \frac{\lambda^2_4 \sqrt{n}}{3(\vartheta+\lambda_4)}\right\}.~~~~~\eeqlb
Next we show $\varlimsup_{n\ra\infty} d_n\mbfP\left(Q^c_n\right)<+\infty.$
Note that $$\mbfP\left(Q^c_n\right)\leq\mbfP\left(H^c_n\right)+\mbfP\left(J^c_n\right)+\mbfP(\tilde{J}^c_n)+\mbfP\left(I^c_n\right).$$
Since the random environment is i.i.d., we see for any $n,$ $\{\rho_{n,i}, i\in\bfN^+\}$ is an i.i.d. sequence, which means that
$$\mbfP\left(\min_{i\leq q^*_n}\rho_{n,i}<q_n\right)\leq q^*_n\mbfP(\rho_{n,1}<q_n).$$
Recalling the definition of $\rho_{n,1}$ we see $\mbfP(\rho_{n,1}<q_n)\leq \mbfP\left(\max_{i\leq q_n}|M_i|> \delta\sqrt{n}\right)$. According to \eqref{sect-2} and the first equality in \eqref{prop1}, we see $\mbfE M_1=0.$
Recall the notation $\sigma^2:=\mathbf{E}\left[\left(\kappa_1(\vartheta)-\vartheta\kappa'_1(\vartheta)\right)^2\right]$ and thus $\sigma^2=\mathbf{E}\left(M^2_1\right).$
Now we take $l_n:=q_n$ and $m_n:=\delta\sqrt{n}$, and then the range of $c_2$ in \eqref{rangec2} and assumption \eqref{c2'a} ensure that $(m_n,l_n)$ is a STA couple at level $\left(\lambda_0, 1+\frac{c_2\delta^2}{(\lambda_0-2)\sigma^2}\right)$ w.r.t. $M_1$ for $n$ large enough.
According to Corollary \ref{sat},  we can find a constant $c_4>0$ such that
 $\mbfP\left(\max_{i\leq q_n}|M_i|> \delta\sqrt{n}\right)\leq c_4q_nn^{-\lambda_0/2}$ and hence
\begin{eqnarray}\label{hnc}
\mbfP(H^c_n)=q^*_n\mbfP(\rho_{n,1}<q_n)\leq c_4q^*_nq_nn^{-\frac{\lambda_0}{2}}\leq 2c_2c_4n^{1-\frac{\lambda_0}{2}}\log n\leq 2 c_2 c_4n^{1-\frac{\lambda_0}{2}+\varepsilon}
\end{eqnarray}
holds for $n$ large enough.

According to Lemma \ref{mto}, \eqref{c2'b} is equivalent to
\begin{eqnarray}\label{c2'bT}
\mbfE\left(\left[\mbfE_{\L}\left(|U_1|^{\lambda_2}\right)\right]^{\lambda_1}\right)<+\infty.
\end{eqnarray}
Hence we see $\mbfE(\Gamma^{\lambda_1\lambda_2/2}_1)<+\infty$ by Jensen's inequality (recalling that $\lambda_2>2$). From Lemma \ref{mto} and the definition of $\sigma^2_*$ we see $\sigma^2_*=\mbfE(\Gamma_1).$ By a direct calculation we see that $(c_3q_n,2n)$ is a STA couple at level $(\lambda_1\lambda_2/2,3)$ w.r.t. $\Gamma_1$ as long as $n$ is large enough. Therefore, Corollary \ref{sat} tells that we can find a constant $c_5>0$ such that
\begin{eqnarray}\label{jnc}\mbfP(J^c_n):=\mbfP\left(\max_{i\leq 2n}|\Gamma_{i}-\sigma^2_*i|> c_3q_n\right)\leq c_5nq^{-\lambda_1\lambda_2/2}_n\leq 2c_5n^{1-\lambda_1}\end{eqnarray} for $n$ large enough.

Recall that $\psi_1:=\mbfE_{\L}\left(|U_1|^{\lambda_2}\right).$ Since $$\mbfP(\tilde{J}^c_n)=\mbfP(\psi_{2n}>3\mbfE(\psi_1)n)=\mbfP(\psi_{2n}-2n\mbfE(\psi_1)>n\mbfE(\psi_1)),$$ 
\eqref{c2'bT} and Corollary \ref{sat} tell that \begin{eqnarray}\label{j1nc}\exists c_6>0,~~~\forall n\in\bfN^+,~~~ \mbfP(\tilde{J}^c_n)\leq c_6n^{1-\lambda_1}.\end{eqnarray}
Note that $\left\{\mbfE_{\L}(N(u)^{1+\lambda_4})\leq e^{\frac{\lambda_4\sqrt{n}}{3}}\right\}=\left\{\max\{\mbfE_{\L}(N(u)^{1+\lambda_4}),1\}\leq e^{\frac{\lambda_4\sqrt{n}}{3}}\right\}$ and recall the assumption
$$\mbfE([\log^+\mbfE_{\L}(N(u)^{1+\lambda_4})]^{\lambda_3})+\mbfE(|\kappa_1(\vartheta)|^{\lambda_3})+\mbfE(|\kappa_1(\vartheta+\lambda_4)|^{\lambda_3})<+\infty.$$ Then by Markov inequality we can find a constant $c_7$ such that
\begin{eqnarray}\label{inc}\forall n\in\bfN^+, ~~~\mathbf{P}(I^c_n)\leq c_7n^{1-\frac{\lambda_3}{2}}.\end{eqnarray}
Since we have chosen $\varepsilon$ to satisfy $\min(\frac{\lambda_0}{2}-\varepsilon,\lambda_1,\frac{\lambda_3}{2})>\frac{\lambda_6}{\lambda_6-p-\varepsilon}$, combining with \eqref{hnc}-\eqref{inc} we get $\varlimsup_{n\ra\infty} d_n\mbfP\left(Q^c_n\right)<+\infty.$

{\it Part 3.}~~Recalling \eqref{23<} we see the only rest to prove is
\begin{eqnarray}\label{EQnc}\varlimsup_{n\ra\infty}\sum^{d_n}_{m=1}\mbfP\left(\left|\log \hat{p}_{\L}(T;n,1,a,b,a)\right|^{p+\varepsilon}\geq m, Q_n\right)<+\infty.\end{eqnarray}

By the definition in \eqref{pt1}, we see \begin{eqnarray}\label{end6}\hat{p}_{\L}(T; n,1,a,b,a)\geq p_{\L}(T;n,1,a,b,a)-\sum_{i=1}^n\mbfP_{\L}(\xi_i>e^{\sqrt{n}}),\end{eqnarray}
where \begin{eqnarray}\label{pt10}p_{\L}(T;n,z,a,b,c):=\inf_{|x|\leq a\sqrt{n}}\mbfP_\L
\left(\forall_{i\leq zn, i\in\bfN}~|T_{i}|\leq b\sqrt{n},|T_{\lfloor zn\rfloor}|\leq c\sqrt{n}\big|T_{0}=x\right).\end{eqnarray}
By Lemma \ref{mto}, we have
\begin{eqnarray}\label{zong}\mbfP_{\L}(\xi_i>e^{\sqrt{n}})=\frac{\mbfE_{\L}\left(\1_{\{N(u)>e^{\sqrt{n}}\}}\sum\limits_{j=1}^{N(u)}e^{-\vartheta\zeta_{j}(u)}\right)}{\mbfE_{\L}\left(\sum\limits_{j=1}^{N(u)}e^{-\vartheta\zeta_{j}(u)}\right)} ,~~|u|=i-1.\end{eqnarray}
Then by the method used in the proof of \cite[Corollary 4.3]{LY3}, we get
\begin{eqnarray}\label{pop}\mbfP_{\L}(\xi_i>e^{\sqrt{n}})\leq e^{-\lambda_4v_1\sqrt{n}}\mbfE_{\L}(N(u)^{1+\lambda_4})^{v_1}e^{(1-v_1)\kappa_i(\vartheta+\lambda_4)-\kappa_i(\vartheta)}, ~~|u|=i-1,\end{eqnarray}
where $\lambda_4$ has been introduced in \eqref{c3a} and $v_1:=\frac{\lambda_4}{\vartheta+\lambda_4}.$ Recalling \eqref{tilIn} we see that
\begin{eqnarray}\label{Ilow}\sum_{i=1}^n\mbfP_{\L}(\xi_i>e^{\sqrt{n}})\leq ne^{-\frac{\lambda_4v_1\sqrt{n}}{3}}, ~~\text{on}~~ I_n.\end{eqnarray}
Now we estimate the lower bound of $p_{\L}(T;n,1,a,b,a).$ Define the shift operator~$\mathfrak{T}$~as $\mathfrak{T}\L:=(\L_2,\L_3,\ldots).$ We usually write
$$\mathfrak{T}_0\L:=\L,~~\mathfrak{T}_k:=\mathfrak{T}^{*k},~~\forall k\in\bfN^+~~\text{and~hence}~~\mathfrak{T}_k\L=(\L_{k+1},\L_{k+2},\ldots).$$
Recall \eqref{rangec2-} and let $a_{n,0}:=a,$ $a_{n,i}:=a_{n,i-1}-\frac{\rho_{n,i}(a-\bar{a})}{n}.$ (We remind that $\mbfP(N_n\geq 1)=1$ because of the definition of $\tau_{n,i}$.)
By Markov property we see 
\begin{eqnarray}\label{end7}&&p_{\L}(T;n,1,a,b,a)\no\\&\geq& \prod_{i=1}^{N_n} p_{\mathfrak{T}_{\tau_{n,i-1}}\L}\left(T;n,\rho_{n,i}/n,a_{n,i-1},b,a_{n,i}\right)\times p_{\mathfrak{T}_{\tau_{n,N_n}}\L}\left(T;n,1-\frac{\tau_{_{n,N_n}}}{n},a_{n,N_n},b,a\right),~~\end{eqnarray}
where we agree that $p_{\cdot}\left(T;\cdot,0,\cdot,\cdot,\cdot\right)=1.$
For any $i\leq N_n,$ by the definitions of $\rho_{n,i}, M_n$ and $U_n$, we see
\begin{eqnarray}\label{pi<}&&p_{\mathfrak{T}_{\tau_{n,i-1}}\L}\left(T;n,\rho_{n,i}/n,a_{n,i-1},b,a_{n,i}\right)\no
\\&=&\inf_{|x|\leq a_{n,i-1}\sqrt{n}}\mbfP_{\L}\left(\begin{split}\forall_{\tau_{n,i-1}\leq k\leq \tau_{n,i}}x+U_k-U_{\tau_{n,i-1}} +M_k-M_{\tau_{n,i-1}}\in[-b\sqrt{n},b\sqrt{n}]\\x+U_{\tau_{n,i}}-U_{\tau_{n,i-1}}+M_{\tau_{n,i}}-M_{\tau_{n,i-1}}\in[-a_{n,i}\sqrt{n},a_{n,i}\sqrt{n}]\end{split}\right),~~~~~~~~~\end{eqnarray}
From Theorem \uppercase\expandafter{\romannumeral1} we can find a standard Brownian motion $B$ under $\mbfP_{\L}$ such that
\begin{eqnarray}\label{3.3.8}\exists c_8>0,~~\forall n\in \bfN,~~\mathcal{Q}_n:=\mbfP_{\L}\left(\max_{i\leq 2n}|U_i-B_{\Gamma_i}|\geq \frac{c_3}{2}\delta\sqrt{n}\right)\leq \frac{c_8\psi_{2n}}{n^{\lambda_2/2}},~~{\rm \mathbf{P}\text{-}a.s.,}\end{eqnarray}
and thus on the event $H_n$,
\begin{eqnarray}\label{pi<+}&&p_{\mathfrak{T}_{\tau_{n,i-1}}\L}\left(T;n,\rho_{n,i}/n,a_{n,i-1},b,a_{n,i}\right)\no
\\&\geq&\inf_{|x|\leq a_{n,i-1}\sqrt{n}}\mbfP_{\L}\left(\begin{split}\forall_{\tau_{n,i-1}\leq k\leq \tau_{n,i}}|x+U_k-U_{\tau_{n,i-1}}| \leq (b-\delta)\sqrt{n}\\|x+U_{\tau_{n,i}}-U_{\tau_{n,i-1}}|\leq (a_{n,i}-\delta)\sqrt{n}\end{split}\right)~~~~~~~~~~~(here~k\in\bfN)\no
\\&\geq&\inf_{|x|\leq a_{n,i-1}\sqrt{n}}\mbfP_{\L}\left(\begin{split}\forall_{\tau_{n,i-1}\leq k\leq \tau_{n,i}}|x+B_{\Gamma_{k}}-B_{\Gamma_{\tau_{n,i-1}}}| \leq (b-\delta')\sqrt{n}\\|x+B_{\Gamma_{\tau_{n,i}}}-B_{\Gamma_{\tau_{n,i-1}}}|\leq (a_{n,i}-\delta')\sqrt{n}\end{split}\right)-\mathcal{Q}_n\no
\\&\geq&\inf_{|x|\leq a_{n,i-1}\sqrt{n}}\mbfP_{\L}\left(\begin{split}\forall_{s\in[0, \eta_{n,i}] }|x+B_{s+\Gamma_{\tau_{n,i-1}}}-B_{\Gamma_{\tau_{n,i-1}}}| \leq (b-\delta')\sqrt{n}\\|x+B_{\Gamma_{\tau_{n,i}}}-B_{\Gamma_{\tau_{n,i-1}}}|\leq (a_{n,i}-\delta')\sqrt{n}\end{split}\right)-\mathcal{Q}_n~~~~(here~s\in\bfR^+)\no
\\&=&\inf_{|x|\leq a_{n,i-1}\sqrt{n}}\mbfP_{\L}\left(\begin{split}\forall_{s\in[0, \eta_{n,i}] }|B_s| \leq (b-\delta')\sqrt{n},|B_{\eta_{n,i}}|\leq (a_{n,i}-\delta')\sqrt{n}\end{split}\big|B_0=x\right)-\mathcal{Q}_n\no
\\&=&\inf_{|x|\leq a_{n,i-1}}\mbfP_{\L}\left(\begin{split}\forall_{s\in[0, \eta_{n,i}/n] }|B_s| \leq b-\delta',|B_{\eta_{n,i}/n}|\leq a_{n,i}-\delta'\end{split}\big|B_0=x\right)-\mathcal{Q}_n,\end{eqnarray}
where \begin{eqnarray}\label{delta'}\delta':=(1+c_3)\delta,~~ \eta_{n,i}:=\Gamma_{\tau_{n,i}}-\Gamma_{\tau_{n,i-1}}.\end{eqnarray}
Recall that $a_{n,i-1}-a_{n,i}=(a-\bar{a})\rho_{n,i}/n$. From the estimate on ``$k(t)$" in \cite[(3.6)-(3.10)]{LY201801} \footnote{Note that $b-\delta'-a_{n,i-1}\geq b-\delta'-a$. Therefore, if the $\delta_1, \delta_2, \epsilon$ in \cite{LY201801} are replaced by $a_{n,i-1}-a_{n,i}+\delta', b-\delta'-a, c_3\delta$ in our context respectively, then the second and third inequalities in \eqref{rangec2-} ensure the truth of \eqref{Plow}.} we see that there exist $c_{9}, c_{10},c_{11}>0$ such that
\begin{eqnarray}\label{Plow}&&\inf_{|x|\leq a_{n,i-1}}\mbfP_\L\left(\forall_{s\leq \frac{\eta_{n,i}}{n}} |B_s|\leq b-\delta', \left|B_{\frac{\eta_{n,i}}{n}}\right|\leq a_{n,i}-\delta'\big|B_0=x\right)\no
\\&\geq&\1_{\{\frac{\eta_{n,i}}{n}\leq c_9\}}\exp\left\{-\frac{(\frac{\rho_{n,i}(a-\bar{a})}{n}+\delta'')^2 n}{2\eta_{n,i}}\right\}+\1_{\{\frac{\eta_{n,i}}{n}> c_9\}}c_{10}\exp\left\{-\frac{c_{11}\eta_{n,i}}{n}\right\},~~~~~\end{eqnarray}
where \begin{eqnarray}\label{delta''}\delta'':=(1+2c_3)\delta.\end{eqnarray} From the definition of $\tau_{n,i}$ we see $\rho_{n,i}\leq n.$ Moreover, on the event $H_n\cap J_n$, we have
\begin{eqnarray}\label{etaup}\forall i\leq N_n+1,~~\eta_{n,i}\leq \sigma_*^2\tau_{n,i}+c_3q_n-(\sigma_*^2\tau_{n,i-1}-c_3q_n)\leq \sigma_*^2n+2c_3q_n\end{eqnarray}
and (recalling the first inequality in \eqref{range+})
\begin{eqnarray}\label{etalow}\eta_{n,i}\geq \sigma^2_*\tau_{n,i}-c_3q_n-(\sigma^2_*\tau_{n,i-1}+c_3q_n)\geq \sigma^2_*\rho_{n,i}-2c_3q_n> \max\left(\frac{\sigma^2_*\rho_{n,i}}{2}, (\sigma^2_*-2c_3)q_n\right).~~~~~~\end{eqnarray}Then on $H_n\cap J_n,$ it is true that
\begin{eqnarray}\frac{(\frac{\rho_{n,i}(a-\bar{a})}{n}+\delta'')^2 n}{2\eta_{n,i}}&=&\frac{n\delta''^2}{2\eta_{n,i}}+\frac{(a-\bar{a})\rho_{n,i}\delta''}{\eta_{n,i}}+\frac{(a-\bar{a})^2\rho^2_{n,i}}{2n\eta_{n,i}}\no
\\&\leq&\frac{n\delta''^2}{2\eta_{n,i}}+\frac{2(a-\bar{a})\delta''}{\sigma^2_*}+\frac{(a-\bar{a})^2}{\sigma^2_*}.\no\end{eqnarray}
According to the above inequality and \eqref{etaup}, there exists $c_{12}$ such that on $H_n\cap J_n$,
\begin{eqnarray}&&\inf_{|x|\leq a_{n,i-1}}\mbfP_{\L}\left(\begin{split}\forall_{s\in[0, \eta_{n,i}/n] }|B_s| \leq b-\delta',|B_{\eta_{n,i}/n}|\leq a_{n,i}-\delta'\end{split}\big|B_0=x\right)\no
\geq 2c_{12}\exp\left\{-\frac{n\delta''^2}{2\eta_{n,i}}\right\}.\end{eqnarray}
Combining the above inequality with \eqref{3.3.8} and \eqref{pi<+}, we see on $H_n\cap J_n \cap \tilde{J}_n,$
\begin{eqnarray}\label{a-b}\forall i\leq N_n,~p_{\mathfrak{T}_{\tau_{n,i-1}}\L}\left(T;n,\rho_{n,i}/n,a_{n,i-1},b,a_{n,i}\right)\geq 2c_{12}\exp\left\{-\frac{n\delta''^2}{2\eta_{n,i}}\right\}-\frac{3c_8\mbfE(\psi_1)n}{n^{\lambda_2/2}}.
\end{eqnarray}
Recalling \eqref{etalow}, the first inequality in \eqref{rangec2} and the definition $q_n:=\frac{n}{c_2\log n}$,
we see 
$$\exp\left\{-\frac{n\delta''^2}{2\eta_{n,i}}\right\}\geq n^{-\frac{\delta''^2c_2}{2(\sigma^2_*-2c_3)}}
~~\text{and}~~\frac{n}{n^{\lambda_2/2}}=o\left(\exp\left\{-\frac{n\delta''^2}{2\eta_{n,i}}\right\}\right).$$
Hence on $H_n\cap J_n \cap \tilde{J}_n$, \eqref{a-b} and \eqref{etalow} tell that for any $i\leq N_n,$
\begin{eqnarray}\label{a-b1} p_{\mathfrak{T}_{\tau_{n,i-1}}\L}\left(T;n,\rho_{n,i}/n,a_{n,i-1},b,a_{n,i}\right)\geq c_{12}\exp\left\{-\frac{n\delta''^2}{2\eta_{n,i}}\right\}\geq c_{12}\exp\left\{-\frac{n\delta''^2}{\sigma^2_*\rho_{n,i}}\right\}~~
\end{eqnarray}
holds for $n$ large enough.

On the other hand, the relationship $\bar{a}<a-\delta'$ in \eqref{rangec2-} and the definition of $N_n$ imply that
\begin{eqnarray}\label{end7+} &&p_{\mathfrak{T}_{\tau_{n,N_n}}\L}\left(T;n,1-\frac{\tau_{_{n,N_n}}}{n},\bar{a},b,a\right)\no
\\&\geq& p_{\mathfrak{T}_{\tau_{n,N_n}}\L}\left(T;n,\frac{\rho_{_{n,N_n+1}}}{n},\bar{a},a,a\right)\no
\\&\geq& p_{\mathfrak{T}_{\tau_{n,N_n}}\L}\left(U;n,\frac{\rho_{_{n,N_n+1}}}{n},\bar{a},a-\delta,a-\delta\right)\no
\\&\geq& \inf_{|x|\leq \bar{a}}\mbfP_\L\left(\forall_{s\leq \frac{\eta_{n,N_n+1}}{n}}~ |B_s|\leq a-\delta'|B_0=x\right)-\mathcal{Q}_n\end{eqnarray}
holds on the event $H_n\cap J_n.$ 
It is plain to see $\tau_{n,N_n+1}\leq 2n$ because of the definitions of $N_n$ and $\tau_{n,i}.$ So we have
$$\eta_{_{n,N_n+1}}=\Gamma_{\tau_{n,N_n+1}}-\Gamma_{\tau_{n,N_n}}\leq \Gamma_{\tau_{n,N_n+1}}\leq 2\sigma_*^2n+c_3q_n,~~{\rm on}~~J_n$$
and thus $$p_{\mathfrak{T}_{\tau_{n,N_n}}\L}\left(T;n,1-\frac{\tau_{_{n,N_n}}}{n},\bar{a},b,a\right)\geq \mbfP_\L\left(\forall_{s\leq 3\sigma^2_*}~ |B_s|\leq a-\delta'-\bar{a}|B_0=0\right)-\mathcal{Q}_n, ~~{\rm on}~~H_n\cap J_n$$
for $n$ large enough. Note that $\mbfP_\L\left(\forall_{s\leq 3\sigma^2_*}~ |B_s|\leq a-\delta'-\bar{a}|B_0=0\right)$ is a positive constant, hence for $n$ large enough, there exists $c_{13}>0$ such that
\begin{eqnarray}\label{pin}p_{\mathfrak{T}_{\tau_{n,N_n}}\L}\left(T;n,1-\frac{\tau_{_{n,N_n}}}{n},\bar{a},b,a\right)\geq c_{13},~~{\rm on}~~H_n\cap J_n\cap \tilde{J}_n.\end{eqnarray}
From \eqref{end7}, \eqref{a-b1} and \eqref{pin}, we have
\begin{eqnarray}&&~~ \label{Ulow}p_{\L}(T;n,1,a,b,a)
\geq c_{13}\prod_{i=1}^{N_n}\left(c_{12}\exp\left\{-\frac{n\delta''^2}{\sigma^2_*\rho_{n,i}}\right\}\right),~~ \text{on}~~H_n\cap J_n\cap \tilde{J}_n\end{eqnarray}
for $n$ large enough. Note that $\1_{H_n}\sum_{i=1}^{N_n}\frac{n}{\rho_{n,i}}\leq N_n\frac{n}{q_n}\leq\frac{n^2}{q^2_n}$ and $\frac{n^2}{q^2_n}=o(\sqrt{n})$.
Hence from \eqref{end6}, \eqref{Ilow} and \eqref{Ulow} we can find a constant $c_{14}>0$ such that
$$\hat{p}_{\L}(T; n,1,a,b,a)\geq \prod_{i=1}^{N_n}\exp\left\{-\frac{c_{14}n}{\rho_{n,i}}\right\},~~\text{on}~~ H_n\cap J_n\cap \tilde{J}_n\cap I_n$$
for $n$ large enough. Finally we completes the proof of Therem \ref{Lp} by the following two propositions.

\begin{prop}\label{lem1}
There exists $\iota\in(0,1)$ such that $\varlimsup_{n\ra +\infty}\mbfP(N_n\geq m)\leq \iota^{m}.$ 
\end{prop}
{\it Proof~of~Proposition \ref{lem1}} By the definition of $\{\tau_{n,i}\}$ and $\{\rho_{n,i}\}$, we see for fixed $n,$ $\rho_{n,1}, \rho_{n,2},...,$ $\rho_{n,i},...$ is an i.i.d. positive random sequence.  Let $F_{(n)}$ be the distribution of $\rho_{n,i}$. Note that
$$\mbfP(N_n\geq m)=\mbfP(\tau_{n,m}\leq n)=F^{*m}_{(n)}(n)\leq [F_{(n)}(n)]^m=[\mbfP(\tau_{n,1}\leq n)]^m.$$
Moreover, according to the classical invariance principle (see \cite{D1951}), we can find a constant $l_1\in(0,1)$ such that
\begin{eqnarray}
\mbfP(\tau_{n,1}>n)&=&\mbfP\left(\max_{i\leq n}|M_i|\leq\delta\sqrt{n}\right)>l_1
\end{eqnarray}
for $n$ large enough, which completes the proof.
\qed
\begin{prop}\label{lem111} 
For any given positive constant $u$, we have
$$\varlimsup_{n\ra +\infty}\mbfE\left(\left(\sum_{i=1}^{N_n}\frac{n}{\rho_{n,i}}\right)^{u}\1_{H_n}\right)<+\infty.$$
\end{prop}
{\it Proof~of~Proposition \ref{lem111}}
Recall that $H_n:=\mbfP\left(\min_{i\leq q^*_n}\rho_{n,i}\geq q_n\right), q_n:=\frac{n}{c_2\log n}$ and $q^*_n:=\left(\lceil\frac{n}{q_n}\rceil\right)^2,$ which means $\1_{H_n}\sum_{i=1}^{N_n}\frac{n}{\rho_{n,i}}\leq q^*_n.$ Then we have
\begin{eqnarray}\label{vs0}\mbfE\left(\left(\sum_{i=1}^{N_n}\frac{n}{\rho_{n,i}}\right)^{u}\1_{H_n}\right)&\leq& 1+\sum_{m=1}^{q^{*^u}_n}\mbfP\left(\left(\sum_{i=1}^{N_n}\frac{n}{\rho_{n,i}}\right)^{u}\1_{H_n}\geq m\right)\no
\\&=&1+\sum_{m=1}^{q^{*^u}_n}\mbfP\left(\left(\sum_{i=1}^{N_n}\frac{n}{\rho_{n,i}}\right)\1_{H_n}\geq m^{\frac{1}{u}}\right).\end{eqnarray}
Choose a constant $\varsigma\in(0,\frac{1}{2})$ and denote $G_{n,m}:=H_n\cap\{N_n\leq \lfloor\varsigma m^{1/u}\rfloor\}.$ It is true that
\begin{eqnarray}\label{vs}\mbfP\left(\left(\sum_{i=1}^{N_n}\frac{n}{\rho_{n,i}}\right)\1_{H_n}\geq m^{\frac{1}{u}}\right)&\leq&\mbfP\left(\sum_{i=1}^{N_n}\frac{n}{\rho_{n,i}}\geq m^{\frac{1}{u}}, G_{n,m}\right)+\mbfP\left(N_n> \lfloor\varsigma m^{\frac{1}{u}}\rfloor\right)\no
\\&\leq& \mbfP\left(\sum_{i=1}^{\lfloor\varsigma m^{1/u}\rfloor}\frac{n}{\rho_{n,i}}\geq m^{\frac{1}{u}}, H_n\right)+\iota^{\lfloor\varsigma m^{1/u}\rfloor}.\end{eqnarray}
Note that $\varsigma\in(0,\frac{1}{2})$ and $m\leq q^{*^u}_n$ mean that $\lfloor\varsigma m^{1/u}\rfloor\leq q^*_n$ for $n$ large enough. Then by Markov inequality we see
\begin{eqnarray}\label{vs1}\forall s>0,~~ \mbfP\left(\sum_{i=1}^{\lfloor\varsigma m^{1/u}\rfloor}\frac{n}{\rho_{n,i}}\geq m^{\frac{1}{u}}, H_n\right)&\leq&e^{-sm^{1/u}}\mbfE\left(e^{s\sum_{i=1}^{\lfloor\varsigma m^{1/u}\rfloor}\frac{n}{\rho_{n,i}}}\1_{H_n}\right)\no
 \\&\leq&e^{-sm^{1/u}}\left[\mbfE\left(e^{\frac{sn}{\rho_{n,1}}}\1_{\{\rho_{n,1}\geq q_n\}}\right)\right]^{\lfloor\varsigma m^{1/u}\rfloor}.\end{eqnarray}
Now we need to estimate the tail of $\frac{n}{\rho_{n,1}}.$ For any $x>1$, by the definition of $\rho_{n,1}$ in \eqref{rho*} we see
$$\mbfP\left(\frac{n}{\rho_{n,1}}>x\right)=\mbfP\left(\frac{n}{x}>\rho_{n,1}\right)\leq \mbfP\left(\max_{i\leq \lceil\frac{n}{x}\rceil}|M_i|\geq \delta\sqrt{n}\right).$$
Recall that $\sigma^2=\mbfE(M^2_1).$ By the method used in the proof of Corollary \ref{sat} (noting that here we use \eqref{citeprop} in the proof of Corollary \ref{sat} rather than the result of Corollary \ref{sat} directly), we can find constants $l_2$ and $l_3$ such that
$$\forall n\in\bfN^+,~\mbfP\left(\max_{i\leq \lceil\frac{n}{x}\rceil}|M_i|>\delta\sqrt{n}\right)\leq l_2e^{-\frac{x\delta^2}{3\sigma^2}}+l_3\lceil\frac{n}{x}\rceil n^{-\frac{\lambda_0}{2}}.$$
Note that  $e^{\frac{sn}{\rho_{n,1}}}\1_{\{\rho_{n,1}\geq q_n\}}\leq n^{sc_2}.$ By the above inequality we see
\begin{eqnarray}\forall s\leq 1,~\forall n\in\bfN^+,~\mbfE\left(e^{\frac{sn}{\rho_{n,1}}}\1_{\{\rho_{n,1}\geq q_n\}}\right)&\leq&3+ \sum^{\lceil n^{sc_2}\rceil}_{k=3}\mbfP\left(e^{\frac{sn}{\rho_{n,1}}}\1_{\{\rho_{n,1}\geq q_n\}}\geq k\right)\no
\\&\leq&3+ \sum^{\lceil n^{sc_2}\rceil}_{k=3}\mbfP\left(\frac{n}{c_2\log n}\leq\rho_{n,1}\leq\frac{sn}{\log k}\right)\no
\\&\leq&3+ \sum^{\lceil n^{sc_2}\rceil}_{k=3}\mbfP\left(\frac{n}{\rho_{n,1}}\geq\frac{\log k}{s}\right)~~~\text{\big(noting~that} ~\frac{\log k}{s}>1\text{\big)}\no
\\&\leq&3+ \sum^{\lceil n^{sc_2}\rceil}_{k=3}\left(l_2k^{-\frac{\delta^2}{3s\sigma^2}}+l_3\lceil\frac{ns}{\log k}\rceil n^{-\frac{\lambda_0}{2}}\right)\no
\\&\leq&3+ \left(l_2\sum^{+\infty}_{k=3}k^{-\frac{\delta^2}{3s\sigma^2}}\right)+\lceil n^{sc_2}\rceil l_3\lceil\frac{ns}{\log 3}\rceil n^{-\frac{\lambda_0}{2}}.\no
\end{eqnarray}
Recall that $\lambda_0>3,$ hence we can choose $s\in\left(0, \min\left\{\frac{\delta^2}{3\sigma^2}, \frac{\lambda_0-2}{2c_2},1\right\}\right)$. Then there exists a constant $l_4>1$ such that
\begin{eqnarray}\label{l4}\sup_n\mbfE\left(e^{\frac{sn}{\rho_{n,1}}}\1_{\{\rho_{n,1}\geq q_n\}}\right)<l_4.\end{eqnarray}
Finally, combining \eqref{l4} with \eqref{vs0}-\eqref{vs1}, we complete the proof by choosing $$\varsigma\in \left(0, \min\left(\frac{s}{\log l_4},\frac{1}{2}\right)\right).$$ ~~~~~~~~~~~~~~~~~~~~~~~~~~~~~~~~~~~~~~~~~~~~~~~~~~~~~~~~~~~~~~~~~~~~\qed



\section{The $L^p ~(p\geq 1)$ convergence of small deviation principle for RWre}
Using a method which is similar to the one used in the proof of Theorem \ref{Lp}, we can obtain the $L^p (p\geq 1)$ convergence of the scaling small deviation of a random walk with random environment in time (RWre).
\subsection{Model and result}
First we quote the definition of RWre in Lv, Hong \cite{Lv201802}. Denote $\mu:=\{\mu_n\}_{n\in\bfN^+}$ an i.i.d. sequence taking values in the space of probability measures on $\bfR.$ Conditioned to a realization 
of $\mu,$ we sample $\{X_n\}_{n\in\bfN^+}$ a sequence of independent random variables such that for every $n,$ the law of $X_n$ is the realization of $\mu_n.$
Set \beqlb\label{sec-sn} S_0=x\in\bfR,~~S_n:=S_0+\sum_{i=1}^{n}X_i.\eeqlb 
We call $\{S_n\}_{n\in\bfN}$ a RWre\footnote{Note that this process is not the
 random walk in random environment which has been well-studied in Zeitouni \cite{Z2004}, where the random environment varies in space but for our model, the random environment varies in time.} (with the random environment $\mu$), which also has two laws to be considered. We write $\mbfP_{\mu}$ (quenched law) for the law of $\{S_n\}_{n\in\bfN}$ conditionally on $\mu$ and $\mathbf{P}$ the joint law of $\{S_n\}_{n\in\bfN}$ and $\mu.$ The marginal distribution of $\{S_n\}_{n\in\bfN}$ under $\mbfP$ is called an annealed law. Slightly abusing notation we also write $\mbfP$ for the annealed law.

We remind that the associated walk $\{T_n\}$ introduced in \eqref{shift0} is exactly a RWre. More precisely, as the definition in the previous paragraph, the random environment for $\{T_n\}$ should be the $\{\tau_{n,\L}\}_{n\in\bfN^+}$ introduced in \eqref{m1} with $A=+\infty$. But we can also say $\{T_n\}$ is a random walk with random environment $\L$ as $\{\tau_{n,\L}\}_{n\in\bfN^+}$ is totally determined by $\L$. From this point of view, we can see the close connection between BRWre and RWre.

 The small deviation problem focuses on the probability that a stochastic process has fluctuations below its natural scale. For the small deviation of RWre, we have obtained the sufficient conditions for convergence in probability / almost surely as follows.

\begin{thm}\label{mog}(Lv $\&$ Hong~\cite[Theorem 2]{Lv201802})\label{sdpas}
Assume that
there exist constants $\beta_0>2, \beta_1>2, \beta_2>2$ such that
\begin{eqnarray}\label{m1u1}\mbfE M_1=0,~\mbfE(U_1^2)>0,~\mbfE(|M_{1}|^{\beta_0})+\mbfE\left(\left[\mbfE_\mu(|U_1|^{\beta_2})\right]^{\beta_1}\right)<+\infty,\end{eqnarray}
where $M_{n}:=\mbfE_\mu(S_n), ~U_{n}:=S_n-\mbfE_\mu(S_n).$
Let $\{t_n\}_{n\in\bfN}$ be a sequence of non-negative integers 
and $\bar{t}_n:= t_n+n.$~Set $b_1<a_1\leq a_2<b_2,~b_1\leq a'_1<a'_2\leq b_2$ and $\{y_n\}$ a positive sequence satisfying $y_n=o(\sqrt{n}).$ 

{\rm (a)} The following convergence\beqlb\label{upper0}\lim\limits_{n\rightarrow +\infty}\frac{\sup\limits_{x\in\bfR}\log \mbfP_\mu
\left(\forall_{t_n\leq i\leq \bar{t}_n}~S_{i}\in[b_1y_n, b_2y_n]|S_{t_n}=x\right)}{ny_n^{-2}}= \frac{-\sigma^2_Q}{(b_2-b_1)^2}\gamma\left(\frac{\sigma_A}{\sigma_Q}\right),\eeqlb
\beqlb\label{lower0}&&
\lim\limits_{n\rightarrow +\infty}\frac{\inf\limits_{x\in[a_1 y_n, a_2 y_n]}\log \mbfP_\mu
\left(\begin{split}\forall_{t_n\leq i\leq \bar{t}_n}~ S_{i}\in[b_1y_n, b_2y_n],\\~S_{\bar{t}_n}\in [a'_1 y_n, a'_2 y_n]\end{split}\Bigg|S_{t_n}=x\right)}{ny_n^{-2}}= \frac{-\sigma^2_Q}{(b_2-b_1)^2}\gamma\left(\frac{\sigma_A}{\sigma_Q}\right)~~\eeqlb
hold in probability $\mathbf{P}$ if $ny^{-\min(\beta_0, 2\beta_1)}_n\ra 0$, where
$\sigma^2_{_A}:=\mbfE(M_1^2),\sigma^2_{_Q}:=\mbfE(U_1^2)$ and $\gamma$ is a real-valued function defined in Lv \cite{LY201801}.

{\rm (b)} The convergence in \eqref{upper0} and \eqref{lower0} hold in the sense of ~${\rm \mathbf{P}\text{-}a.s.}$ when $\sum_{n=1}^{+\infty}ny^{-\min(\beta_0, 2\beta_1)}_n< +\infty.$
\end{thm}

The following theorem is the main result in this section.
\begin{thm}\label{sdplp}
Assume that \eqref{m1u1} holds and $\frac{\sigma^2_A}{\sigma^2_Q}< \frac{\beta_2-2}{\beta_0-2}$. If there exist $p\geq 1, t>0, \beta_4>\beta_3>0, \beta_5\geq 1$ such that
\begin{eqnarray}&&\label{lpf1}\mbfE\left(|\log\mbfP_\mu(S_{1}\in [-\beta_4, -\1_{a_2\geq a'_2}\beta_3]|S_{0}=0)|^{\beta_5}\right)\no
\\&&~~~~~~~~~~+\mbfE\left(|\log\mbfP_\mu(S_{1}\in [\1_{a_1\leq a'_1}\beta_3,\beta_4]|S_{0}=0)|^{\beta_5}\right)<+\infty\end{eqnarray} and
\begin{eqnarray}\label{lpf2}t\beta_5\geq 2,~~ p\in\left[\beta_5-\frac{1}{t},\beta_5\right),~~\min\left(\frac{\beta_0}{2}, \beta_1\right)>\frac{\beta_5}{\beta_5-p},\end{eqnarray} then for any $\{y_n\}$ satisfying $ny^{-\min(\beta_0, 2\beta_1)}_n\ra 0$, the convergence in \eqref{upper0} and \eqref{lower0} hold in $L^p.$
\end{thm}
\begin{rem}\label{rema3}
Observing the assumptions in Theorem \ref{sdpas}, for a constant $x$ near $0$, we see that it doesn't matter to the almost surely convergence of the small deviation whether $\mbfP_{\mu}(S_1\leq x)$ is too large or too small. However, we add some restrictions on the tail of $\mbfP_{\mu}(S_1\leq x)$ in \eqref{lpf1} when we consider the $L^p$ convergence of the small deviation. In fact, by an argument similar to Remark \ref{rema2}, we can also see that \begin{eqnarray}\label{r33}\exists x>0
,~~\mbfP(\mbfP_{\mu}(S_1\geq x)=1)=\mbfP(\mbfP_{\mu}(S_1\leq -x)=1)=0\end{eqnarray} is a necessary condition for Theorem \ref{sdplp}. 
From this point we can also see the difference between $L^p$ convergence and almost surely convergence of the quenched small deviation probability of RWre.
\end{rem}

  According to Lemma 
  \ref{mto}, we see the close connection between BRWre and RWre
  and hence many questions on BRWre can also be translated into the corresponding questions on RWre. For example, the random ballot theorem for RWre (\cite[Theorem 1.11]{MM2015}) is a key tool to obtain the asymptotic behavior of the maximal displacement of BRWre, see \cite{MM2016}; the small deviation theorem for RWre (Theorem \ref{sdpas}) plays a key role during we prove some properties of survival probability of BRWre with a barrier, see \cite{LY3}\footnote{Actually, \eqref{sect-2}, \eqref{c2'a} and \eqref{c2'b} are set to ensure that the associated RWre 
  can satisfy the assumptions in Theorem \ref{sdpas}.}. In fact, the estimates on various small deviation probabilities have been used in the study of BRW and some generalizations of BRW more than once, e.g. the BRW with barrier (see \cite{AJ2011,BJ2012,GHS2011}), the $N-$BRW (see \cite{M2017}) and the maximal displacement of a branching random walk in a vary (but non-random) environment (see \cite{M2015b}). Therefore, we believe that there will be potential applications of Theorem \ref{sdplp} in the study of BRWre. Moreover, the $L^p$ convergence for the quenched probability of random model in random environment is of independent interest.

  In view of the applications (mentioned in the previous paragraph) of small deviation to the study of BRWre, we tend to emphasize the importance of \eqref{goal}. According to the Markov property of the RWre, we could view the probability in \eqref{goal} as the ``atom" of the small deviation probability, or more mathematically speaking, the small deviation probabilities in \eqref{upper0}, \eqref{lower0} and \eqref{Bn} can be bounded from above and below by the product of a series of probabilities like the one in \eqref{goal} or \eqref{pt10}. Therefore, \eqref{goal} is highly expected to be applied in more barrier problems or some other topics in BRWre. At least, \eqref{goal} has been used in one of my current researches---the asymptotic behavior of $\mbfP_{\L}(\mathcal{S})$ as $\varepsilon \downarrow 0$ when the barrier function is $i\varepsilon-\vartheta^{-1} K_i$.\footnote{This problem stems from my previous work
  \cite[Theorem 2.5]{LY3}, which states that $\mbfP_{\L}(\mathcal{S})>0, {\rm \mathbf{P}-a.s.}$ (resp. $\mbfP_{\L}(\mathcal{S})=0, {\rm \mathbf{P}-a.s.}$) for $\varepsilon>0$ (resp. $\varepsilon=0$). Such problem in the context of the time-homogeneous case (BRW) was studied in Gantert et al. \cite{GHS2011}.}
  On the other hand, the idea (the constructions of $\tau_{n,i}$ and $Q_n$) used in the proof of \eqref{goal} may also work in the study of limit behavior of BRWre (especially in the sense of moment or $L^p$). As the end of this subsection, let us restate \eqref{goal} as the following proposition, in which we give the sufficient conditions for \eqref{goal} through a careful review of Section 3.2.
  \begin{prop}
  Assume that \eqref{sect-2}, \eqref{c2'a}-\eqref{c3a}, \eqref{T>} hold with constants
\begin{eqnarray}\label{assu17}\lambda_0>2,~~\lambda_1\geq2,~~\lambda_{2}>2,~~\lambda_{3}>3,~~\lambda_{4}>0,~~\lambda_{5}\leq-1,~~ \lambda_6\geq 1,~~\frac{\sigma^2}{\sigma^2_*}< \frac{\lambda_2-2}{\lambda_0-2},\end{eqnarray}
and
\begin{eqnarray}\label{T>7}
\mathbf{E}\left(\left[\log^-\mbfE_{\L}\left(\1_{N(\phi)\leq |\lambda_5|}\sum_{i=1}^{N(\phi)}\1_{\{\vartheta \zeta_i(\phi)+ \kappa_1(\vartheta)\in [\lambda_5,0]\}} \right)\right]^{\lambda_6}\right)<+\infty,\end{eqnarray}
If constants $p, t$ satisfy that
 \begin{eqnarray}\label{assu27}t\geq 1,~~t\lambda_6\geq 2,~~ p\in\left[\lambda_6-\frac{1}{t},\lambda_6\right),~~\min\left\{\frac{\lambda_0}{2},\lambda_1,\frac{\lambda_3}{2}\right\}>\frac{\lambda_6}{\lambda_6-p},\end{eqnarray}
then \eqref{goal} holds for some $\varepsilon>0.$
\end{prop}

\subsection{\bf Proof of Theorem \ref{sdplp}}
Since there are many similarities between the proof of Theorem \ref{sdplp} and the proof of Theorem \ref{Lp}, we emphasize mainly the distinctive points in the proof of Theorem \ref{sdplp}. 
First, by an argument similar to \eqref{An} and \eqref{Bn}-\eqref{B+n} \footnote{We remind that the analog of \eqref{mog.2} in the context of this section will not occur the term like $``\prod Z_m"$ but only the term like $``\prod Z_{j,n}"$ because of the setting $b_1<a_1$ and $a_2<b_2$. Hence we redefine $R_{j,n}$ in \eqref{mog.2} by $R_{j,n}:=\lfloor y^2_n\rfloor j.$}, we can see the heart of the proof is also to estimate the upper limit of
\begin{eqnarray}\label{pp+}\mbfE\left(\left|\inf\limits_{x\in[a_1 \sqrt{n}, a_2 \sqrt{n}]}\log \mbfP_\mu
\left(\begin{split}\forall_{i\leq n}~ S_{i}\in[b_1\sqrt{n}, b_2\sqrt{n}],~\\ S_{n}\in [\max(a_1,a'_1)\sqrt{n},\min(a_2,a'_2)\sqrt{n}]\end{split}\Big|S_{0}=x\right)\right|^{\bar{p}}\right)\end{eqnarray}
for some $\bar{p}>p.$ Without loss of generality, in the rest of the proof we only consider the case
$a_1+a_2=b_1+b_2=a_1'+a_2'=0$ and denote $a:=a_2, b:=b_2, a':=a'_2$. Hence we need to show
\begin{eqnarray}\label{pp'}\exists \bar{p}>p, ~\varlimsup_{n\ra+\infty}\mbfE\left(\left|\inf\limits_{|x|\leq a\sqrt{n}}\log \mbfP_\mu
\left(\forall_{i\leq n}~ |S_{i}|\leq b\sqrt{n},~|S_{n}|\leq a'\sqrt{n}\big|S_{0}=x\right)\right|^{\bar{p}}\right)<+\infty~~~\end{eqnarray}
when $a\geq a'.$ In this section we mainly pay attention to the case $a>a'$ since we have experienced the case $a=a'$ in the proof of Theorem \ref{Lp}.




Hereafter, unless specifically mentioned otherwise, we will continue to use the notation in Section 3.2.

We start the discussion under an extra assumption $a-a'<b-a,$ which means that we can find constants $\delta, c_3, a''$ such that $a''\in(0,a'-\delta')$ and $a-a''<b-a.$ (In fact, we can further choose $a''$ according to the way we choose $\bar{a}$ and then in this proof, $a''$ plays the same role as $\bar{a}$ in Section 3.2.)
Let
 $a_{n,0}:=a,$ $a_{n,i}:=a_{n,i-1}-\frac{\rho_{n,i}(a-a'')}{n}.$ Following the notation in \eqref{pt10}, by Markov property we see 
\begin{eqnarray}\label{end7new}p_{\mu}(S;n,1,a,b,a')
&\geq& \prod_{i=1}^{N_n} p_{\mathfrak{T}_{\tau_{n,i-1}}\mu}\left(S;n,\frac{\rho_{n,i}}{n},a_{n,i-1},b,a_{n,i}\right)\no
\\&\times& p_{\mathfrak{T}_{\tau_{n,N_n}}\mu}\left(S;n,1-\frac{\tau_{n,N_n}}{n},a_{n,N_n},b,a'\right).\end{eqnarray}
For the term $p_{\mathfrak{T}_{\tau_{n,i-1}}\mu}\left(S;n,\frac{\rho_{n,i}}{n},a_{n,i-1},b,a_{n,i}\right),$ 
by an argument similar to \eqref{pi<}-\eqref{a-b} we get a result like \eqref{a-b1}. What we should take careful is the term
$p_{\mathfrak{T}_{\tau_{n,N_n}}\mu}\left(S;n,1-\frac{\tau_{n,N_n}}{n},a_{n,N_n},b,a'\right)$. Since $\tau_{n,N_n}$ is random and $a>a'$, it is not sure whether $a_{n,N_n}<a'-\delta'$ or not. If $a_{n,N_n}<a'-\delta'$, then we can follow the discussion in \eqref{end7+}-\eqref{pin} and finally get \eqref{Ulow}. If $a_{n,N_n}\geq a'-\delta'$, then $\tau_{n,N_n}\leq \frac{a-(a'-\delta')}{(a-a'')}n.$ Therefore, $a_{n,N_n}\geq a'-\delta'$ implies that $n-\tau_{n,N_n}\in[\frac{a'-\delta'-a''}{a-a''}n,n]$ (recalling that $\frac{a'-\delta'-a''}{a-a''}>0$), which means that on the event $H_n\cap J_n$, the time span $\nu_n:=\frac{1}{n}(\Gamma_n-\Gamma_{\tau_{n,N_n}})$ in the probability
$$p(B):=\inf_{|x|\leq a_{n,N_n}}\mbfP_{\L}\left(\begin{split}\forall_{s\in[0, \nu_n] }|B_s| \leq b-\delta',|B_{\nu_n}|\leq a'-\delta'\end{split}\big|B_0=x\right)$$
has two positive constants as its lower bound and upper bound. So conditionally on $H_n\cap J_n$, there is a positive constant as the lower bound of $p(B).$ By an argument similar to \eqref{pi<}-\eqref{a-b} we see that on the event $H_n\cap J_n\cap\tilde{J}_n$, \begin{eqnarray}\label{saklambda--}p(S):=p_{\mathfrak{T}_{\tau_{n,N_n}}\mu}\left(S;n,1-\frac{\tau_{n,N_n}}{n},a_{n,N_n},b,a'\right)
&\geq& p(B)-\frac{3c_8\mbfE(\psi_1)n}{n^{\lambda_2/2}}.~~~~~\end{eqnarray}
So conditionally on $H_n\cap J_n\cap \tilde{J}_n$, there is a positive constant as the lower limit of $p(S)$ as long as $a-a'<b-a.$ Then by Propositions 4.1 and 4.2 we finish the proof for the case $a-a'<b-a$.

At last, if $a,b,a'$ do not satisfy $a-a'<b-a,$ we can divide the target probability in \eqref{pp'} to several segments by Markov property. Just note that for any $p\geq 1, \varrho\in(0,b),$ we have
 \begin{eqnarray}&&\mbfE\left(\left|\log\inf\limits_{|x|\leq a \sqrt{n}} \mbfP_\mu
\left(\forall_{i\leq 2n}~ |S_{i}|\leq b\sqrt{n},~|S_{n}|\leq a'\sqrt{n}\big|S_{0}=x\right)\right|^{p}\right)\no
\\&\leq&2^{p-1}\mbfE\left(\left|\log\inf\limits_{|x|\leq a \sqrt{n}} \mbfP_\mu
\left(\forall_{i\leq n}~ |S_{i}|\leq b\sqrt{n},~|S_{n}|\leq \varrho \sqrt{n}\big|S_{0}=x\right)\right|^{p}\right)\no
\\&+&2^{p-1}\mbfE\left(\left|\log\inf\limits_{|x|\leq \varrho \sqrt{n}} \mbfP_\mu
\left(\forall_{i\leq n}~ |S_{i}|\leq b\sqrt{n},~|S_{n}|\leq a'\sqrt{n}\big|S_{0}=x\right)\right|^{p}\right).\no
\end{eqnarray}
Hence we can choose $\varrho_1> \varrho_2> ...> \varrho_k$ such that $a-\varrho_1<b-a, \varrho_k-a'<b-\varrho_k$ and $\varrho_{i-1}-\varrho_i< b-\varrho_{i-1},\forall 2\leq i\leq k$ and then the situation will go back to the case $a-a'<b-a,$ which completes the proof.\qed

Let us add some comments about the case $a<a'$ for \eqref{pp'}. In fact, if $a<a'$, then \eqref{pp'} requires fewer assumptions and the proof of \eqref{pp'} will be easier.~
That is because under this case, we do not need to worry the situation that $\rho_{n,i}$ may be too small to make the last line in \eqref{pi<+} larger than $0$ and thus we do not need to introduce $H_n$. However, considering the case $a<a'$ is meaningless to the small deviation probability in \eqref{lower0}. The reason is as follows. Recall \eqref{lower0} and note that $|a_2 y_n-a'_2y_n|/(ny^{-2}_n)$ must be $o(y_n)$ because of the assumption $y_n=o(\sqrt{n}).$ (We remind that the assumption $y_n=o(\sqrt{n})$ is absolutely necessary because what we consider is the samll deviation principle but not the moderate deviation or large deviation.)
Therefore, even though $a_2<a_2'$ in \eqref{lower0}, after a decomposition similar to \eqref{end7new}, we also need to prove
\eqref{pp'} under the condition $a'y_n\leq ay_n+o(y_n),$ which is almost the same as the case $a\geq a'.$
 \ack
The author is greatly indebt to Professor Wenming Hong for his helpful advice and encouragement. The author also would like to thank the referees greatly for the coming careful review and valuable suggestions.
This work is supported by the Fundamental Research Funds for the Central Universities (NO.2232021D-30) and the National Natural Science Foundation of China (NO.11971062).


\section*{Declarations}
{\noindent{\bf Conflict of Interest Statement} The author declares that I have no conflict of interest.
\vspace{0.3cm}

\noindent{\bf Availability of Data and Material (data transparency)} Data and material sharing are not applicable to this article as no data sets and material were generated or analyzed during the current (theoretical) study.
\vspace{0.3cm}

\noindent{\bf Data Availability Statement} No data were generated or analyzed in this article.

\vspace{0.3cm}

\noindent{\bf Code Availability (software application or custom code)} Code sharing is not applicable to this article as no
code was generated during the current (theoretical) study.
}

\end{document}